\documentclass{article}
\usepackage{amsmath, amssymb, amsthm}
\usepackage{epsfig}

\newcommand{\Bt}{\overset{\sim}{B}}
\title{The Freedman group: \\a physical interpretation for the $SU(3)$-subgroup $D(18,1,1;2,1,1)$ of order $648$}
\author{Claire Levaillant}

\begin{document}
\maketitle

\noindent Abstract. We study a subgroup $Fr(162\times 4)$ of $SU(3)$ of order $648$ which is an extension of $D(9,1,1;2,1,1)$ and whose generators arise from anyonic systems. We show that this group is isomorphic to a semi-direct product $(\mathbb{Z}/18\mathbb{Z}\times\mathbb{Z}/6\mathbb{Z})\rtimes S_3$ with respect to conjugation and we give a presentation of the group. We show that the group $D(18,1,1;2,1,1)$ from the series $(D)$ in the existing classification for finite $SU(3)$-subgroups is also isomorphic to a semi-direct product $(\mathbb{Z}/18\mathbb{Z}\times\mathbb{Z}/6\mathbb{Z})\rtimes S_3$, also with respect to conjugation. We show that the two groups $Fr(162\times 4)$ and $D(18,1,1;2,1,1)$ are isomorphic and we provide an isomorphism between both groups. We prove that $Fr(162\times 4)$ is not isomorphic to the exceptional $SU(3)$ subgroup $\Sigma(216\times 3)$ of the same order $648$. We further prove that the only $SU(3)$ finite subgroups from the $1916$ classification by Blichfeldt or its extended version which $Fr(162\times 4)$ may be isomorphic to belong to the $(D)$-series. Finally, we show that $Fr(162\times 4)$ and $D(18,1,1;2,1,1)$ are both conjugate under an orthogonal matrix which we provide.  %more generally that these are the only $SU(3)$ finite subgroups from the classification which have that given structure after we exclude the finite subgroups of $U(2)$.

% contains two more groups, namely two direct products of $\mathbb{Z}_3$ by an exceptional group. We thus add a new finite subgroup of $SU(3)$ to the incomplete classification, moving one more step forward towards a complete classification of all the finite subgroups of $SU(3)$. We foresee that this result of interest to both the mathematical and the physical community will have groundbreaking applications in quantum computation or particle physics.

\section{Introduction and Definition of the group}
Finite subgroups of $SU(3)$ are used in particle physics and topological quantum computation among other fields. The first tentative classification of the finite subgroups of $SU(3)$ appears in the works of G.A. Miller, H.F. Blichfeldt and L.E. Dickson \cite{MB}, dating from $1916$. The part written by Blichfeldt describes in terms of generators all the finite $SU(3)$ subgroups currently known, excluding two more recent ones, namely direct products with a cyclic group of order $3$ of the two smallest non-abelian simple groups, that is $A_5\times\mathbb{Z}/3\mathbb{Z}$ and $PSL(2,7)\times\mathbb{Z}/3\mathbb{Z}$. Since then,
many efforts have been pursued and progress is still made. For instance, only recently in \cite{PO2}, P.O. Ludl determines the structure of the series $(C)$ and $(D)$ and gives in particular a simple example for a (C)-group, namely $C(9,1,1)$, which is neither of the form $\Delta(3n^2)$ nor
of the form $T_n$, thus showing that (C) contains some hitherto unclassified subgroups of
$SU(3)$ (the sub-series $T_n$ and $\Delta(3n^2)$ of $(C)$ were studied in \cite{BO1},\cite{BO2},\cite{FF1},\cite{LU}, while the full series $(C)$ remained unstudied until $2011$ in \cite{PO2}). For a clear historical summary and table of the existing classification, see
\cite{PO2}. The group we study in this paper is not one of the $6$ exceptional
subgroups of $SU(3)$, although it has the same order as $\Sigma(216\times 3)$, that is $648$, and is like $\Sigma(216\times 3)$ a group extension of $D(9,1,1;2,1,1)$ by \cite{LL}. It is rather the group $D(18,1,1;2,1,1)$ from the series $(D)$. Our group arises up to phase as the result of unitary quantum gates obtained by braiding $4$ anyons of topological charge $2$ in the Jones-Kauffman version of $SU(2)$ Chern-Simons theory at level $4$ on the one hand and by fusing
a topological charge $4$ out of the vacuum on the other hand, like shown on the figures below. \\
\begin{center}
\epsfig{file=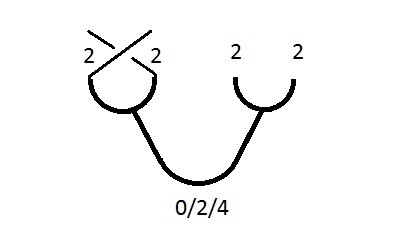, height=2.5cm}\epsfig{file=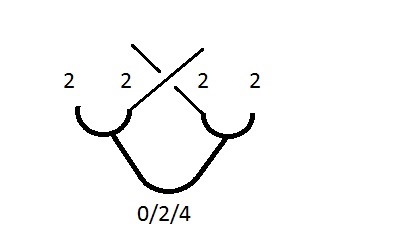, height=2.5cm}\end{center}
\begin{center}\epsfig{file=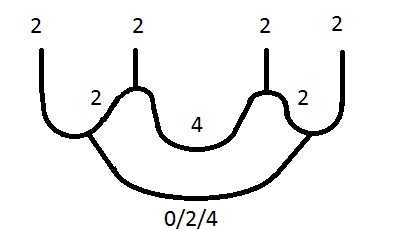, height=2.5cm}\\$\begin{array}{l}\\\textit{The Freedman fusion operation}\end{array}$\end{center}

\noindent The latter operation is due to Mike Freedman \footnote{The Freedman fusion operation was communicated by Mike Freedman to the author during a Skype conversation between Station Q, Santa Barbara and HRI, Allahabad, where the author was visiting. She thanks Harish Chandra Research Institute for being a family.}
and has the
effect of exchanging particles of respective topological charges $0$ and $4$ on the quantum trit, thus allowing to get
an extra quantum gate, namely a permutation matrix up to phase, see Appendix. This extra quantum gate discovered by Mike Freedman provides the third generator of our
group, which we name for that reason $Fr(162\times 4)$. In \cite{GL}, P. O. Ludl studies
the principal series of some exceptional groups, thus establishing their structure. His work shows in particular that $\Sigma(216\times 3)$ is isomorphic to
$$\Delta(54)\rtimes A_4=((\mathbb{Z}/3\mathbb{Z}\times\mathbb{Z}/3\mathbb{Z})\rtimes\,S_3)\rtimes (V_4\rtimes\mathbb{Z}/3\mathbb{Z})$$ As part of our paper, we show that $Fr(162\times 4)$ is isomorphic to
$$(\mathbb{Z}/9\mathbb{Z}\times\mathbb{Z}/3\mathbb{Z}\times V_4)\rtimes S_3$$%The main difference with our paper is that he needs
%an extra generator to define the group and the structure which he gets is less straightforward of use.
Let us now proceed to the definition of $Fr(162\times 4)$ by providing its generators.
Our generators are the following: the two braid matrices from \cite{BL} which we recall below and the fusion matrix, which we will denote by $FUM$.
$$\begin{array}{ccc}G_1=\begin{pmatrix}e^{\frac{7i\pi}{9}}&0&0\\&&\\0&-e^{\frac{4i\pi}{9}}&0\\&&\\0&0&-e^{\frac{7i\pi}{9}}
\end{pmatrix}&G_2=\begin{pmatrix}
-\frac{1}{2}\,e^{\frac{4i\pi}{9}}&\frac{e^{\frac{7i\pi}{9}}}{\sqrt{2}}&\frac{1}{2}\,e^{\frac{4i\pi}{9}}\\
&&\\
\frac{e^{\frac{7i\pi}{9}}}{\sqrt{2}}&0&\frac{e^{\frac{7i\pi}{9}}}{\sqrt{2}}\\
&&\\
\frac{1}{2}\,e^{\frac{4i\pi}{9}}&\frac{e^{\frac{7i\pi}{9}}}{\sqrt{2}}&-\frac{1}{2}\,e^{\frac{4i\pi}{9}}
\end{pmatrix}\end{array}$$$$FUM=\begin{pmatrix}0&0&-e^{\frac{2i\pi}{3}}\\&&\\0&-e^{\frac{2i\pi}{3}}&0\\&&\\-e^{\frac{2i\pi}{3}}&0&0
\end{pmatrix}$$

Our result is the following. For clarity, we denote the group $\mathcal{G}$ of \cite{BL} by $Fr(162)$.
\newtheorem{Theorem}{Theorem}

\begin{Theorem}
The subgroup $Fr(162\times 4)$ of $SU(3)$ generated by the three matrices above is a finite group extension of $D(9,1,1;2,1,1)$ of order $648$ and is isomorphic to $((\mathbb{Z}/3\mathbb{Z}\times\mathbb{Z}/9\mathbb{Z})\times V_4)\ltimes S_3$, where the semi-direct product is with respect to conjugation (the detailed action is provided in the presentation below).  Any element of $Fr(162\times 4)$ can be written uniquely as a product
$A^i\,B^j\,T\,H$ with $A$ and $B$ defined like further below, $0\leq i\leq 8$, $0\leq j\leq 2$, $H$ the identity matrix or one of the matrices $H_i$
with $i\in\lbrace 0,1,2,3,4\rbrace$ whose definitions also appear further below and $T$ an element of the Klein group
$$\mathcal{V}\,=\,\left\lbrace\begin{array}{l} I_3, (FUM)^3=\begin{pmatrix} &&-1\\&-1&\\-1&&\end{pmatrix}, G_1\,(FUM)^3\,G_1^{-1}=\begin{pmatrix} &&1\\&-1&\\1&&\end{pmatrix},\\G_1(FUM)^3\,G_1^{-1}(FUM)^3=\begin{pmatrix} -1&&\\&1&\\&&-1\end{pmatrix}\end{array}\right\rbrace$$
Set $C_{18}=A(FUM)^3$ and $C_6=B\,G_1(FUM)^3\,G_1^{-1}$.\\\\ A presentation of $Fr(162\times 4)$ is as follows.
$$\left\langle\begin{array}{ccc}\begin{array}{l}\\ C_6,C_{18},H_1,H_3\\\end{array}&|&\begin{array}{l}C_6^6=C_{18}^{18}=[C_6,C_{18}]=H_1^2=H_3^3=(H_1H_3)^2=I\\\\
\begin{array}{cc}\begin{array}{l}H_1\,C_6\,H_1^{-1}=C_6^{-1}\,C_{18}^6\\H_1\,C_{18}\,H_1^{-1}=C_6^3\,C_{18}\end{array}&\begin{array}{l}H_3\,C_6\,H_3^{-1}=
C_6\,C_{18}^{-3}\\H_3\,C_{18}\,H_3^{-1}=
C_6\,C_{18}^{10}\end{array}\end{array}\end{array}\end{array}\right\rangle$$
$Fr(162\times 4)$ contains exactly four $3$-Sylow subgroups, one of which is the unique $3$-Sylow subgroup of $Fr(162)$ given by
$$\mathcal{S}_3(F)=\mathcal{N}\,\bigsqcup\,\mathcal{N}H_3\,\bigsqcup\,\mathcal{N}H_3^2$$
where $\mathcal{N}=<A,B>$. The other three are the $\mathcal{V}$-conjugates of $\mathcal{S}_3(F)$. \\
The group $Fr(162\times 4)$ is conjugate to the $SU(3)$ finite subgroup $D(18,1,1;2,1,1)$, under an orthogonal matrix as follows.
$$\begin{pmatrix} 1/\sqrt{2}&0&1/\sqrt{2}\\0&1&0\\-1/\sqrt{2}&0&1/\sqrt{2}\end{pmatrix}Fr(162\times 4)\begin{pmatrix}
1/\sqrt{2}&0&-1/\sqrt{2}\\0&1&0\\1/\sqrt{2}&0&1/\sqrt{2}\end{pmatrix}=D(18,1,1;2,1,1)\\$$
\end{Theorem}
$$\begin{array}{l}\end{array}$$
The paper is organized as follows.
In the next part of the paper, we proceed to the study of the structure of the Freedman group, based on the knowledge
we already have of the structure of its subgroup $D(9,1,1;2,1,1)$, generated by the two braid matrices, see \cite{BL}.
We then study the structure of another $SU(3)$ finite subgroup, namely $D(18,1,1;2,1,1)$. Along the way, we prove in particular that the semi-direct product $\mathbb{Z}_{18}\times\mathbb{Z}_6\rtimes S_3$ arises as a finite subgroup of $SU(3)$. We next show that $Fr(162\times 4)$ is isomorphic to $D(18,1,1;2,1,1)$ by exhibiting an isomorphism between both groups. We further show that the group $Fr(162\times 4)$ is isomorphic to at least four other $D$-groups and that the five $D$-groups are identical. Finally, we show that the two groups $Fr(162\times 4)$ and $D(18,1,1;2,1,1)$ are conjugate. Towards the end of the paper, we explain how we could verify some of our results by using the computer program GAP. %Along the way, we prove the a priori non-trivial result that up to isomorphism, a unique semi-direct product $\mathbb{Z}_6\times\mathbb{Z}_{18}\rtimes S_3$ arises in the series $(D)$.
\section{Structure of the group}

First, we recall some facts from \cite{BL} about the structure of $$D(9,1,1;2,1,1)=\,<G_1,G_2>$$
In \cite{BL} we show that this group is isomorphic to the semidirect product $\mathbb{Z}_9\times\mathbb{Z}_3\rtimes S_3$ with respect to conjugation. The elements of the abelian group $\mathbb{Z}_9\times\mathbb{Z}_3$ are the matrices
$A^i\,B^j$ with $0\leq i\leq 8$ and $0\leq j\leq 3$. The matrices $A$ and $B$ commute and are provided below.

$$A=\begin{pmatrix} \frac{1}{2}\,e^{\frac{5i\pi}{9}}(-1+e^{\frac{i\pi}{3}})&0&\frac{1}{2}\,e^{\frac{5i\pi}{9}}
(1+e^{\frac{i\pi}{3}})\\&&\\0&-e^{\frac{5i\pi}{9}}&0\\&&\\
\frac{1}{2}\,e^{\frac{5i\pi}{9}}(1+e^{\frac{i\pi}{3}})&0&\frac{1}{2}\,e^{\frac{5i\pi}{9}}(-1+e^{\frac{i\pi}{3}})
\end{pmatrix}\\$$

$$B=\begin{pmatrix} \frac{1}{2}\,(1+e^{\frac{2i\pi}{3}})&0&\frac{1}{2}(-1+e^{\frac{2i\pi}{3}})\\
&&\\
0&-e^{\frac{i\pi}{3}}&0\\&&\\
\frac{1}{2}(-1+e^{\frac{2i\pi}{3}})&0&\frac{1}{2}(1+e^{\frac{2i\pi}{3}})\end{pmatrix}$$
In terms of generators, these matrices are

\begin{eqnarray}
A&=&G_1\,G_2^2\,G_1^{-1}\\
B&=&G_1\,G_2^{-2}\,G_1
\end{eqnarray}
It is shown in \cite{BL} that the subgroup of $<G_1,G_2>$ generated by the matrices $A$ and $B$ is a normal
subgroup $\mathcal{N}$ of $D(9,1,1;2,1,1)$. Moreover, there is in $D(9,1,1;2,1,1)$ a subgroup $\mathcal{H}$ isomorphic to the symmetric group $S_3$, whose six matrices are given by the identity matrix, the matrix $$H_0=\begin{pmatrix} -1 &0&0\\0&-1&0\\0&0&1\end{pmatrix}$$
and the non-diagonal matrices $H_i$, $1\leq i\leq 4$ with

\begin{center}$$H_1=\begin{pmatrix} -\frac{1}{2}& -\frac{1}{\sqrt{2}}&-\frac{1}{2}\\&&\\
-\frac{1}{\sqrt{2}}&0&\frac{1}{\sqrt{2}}\\&&\\
-\frac{1}{2}&\frac{1}{\sqrt{2}}&-\frac{1}{2}\end{pmatrix},\qquad\qquad H_2= \begin{pmatrix} -\frac{1}{2}& -\frac{1}{\sqrt{2}}&\frac{1}{2}\\&&\\
-\frac{1}{\sqrt{2}}&0&-\frac{1}{\sqrt{2}}\\&&\\
\frac{1}{2}&-\frac{1}{\sqrt{2}}&-\frac{1}{2}\end{pmatrix},$$\end{center}\begin{center} $$H_3=\begin{pmatrix}\frac{1}{2}& \frac{1}{\sqrt{2}}&-\frac{1}{2}\\&&\\
\frac{1}{\sqrt{2}}&0&\frac{1}{\sqrt{2}}\\&&\\
\frac{1}{2}&-\frac{1}{\sqrt{2}}&-\frac{1}{2}\end{pmatrix},\qquad\qquad H_4=\begin{pmatrix}
\frac{1}{2}& \frac{1}{\sqrt{2}}&\frac{1}{2}\\&&\\
\frac{1}{\sqrt{2}}&0&-\frac{1}{\sqrt{2}}\\&&\\
-\frac{1}{2}&\frac{1}{\sqrt{2}}&-\frac{1}{2}\end{pmatrix}$$\end{center}

\noindent The matrices $H_0$, $H_1$ and $H_2$ are the three elements of order $2$ in the group and the matrices $H_3$ and $H_4$ the two elements of order $3$.\\
We recall from \cite{BL} that the intersection $\mathcal{N}\cap\mathcal{H}$ is trivial, and that
$$<G_1,G_2>=\mathcal{N}.\mathcal{H}$$
And so any element of $D(9,1,1;2,1,1)$ can be uniquely written as the product of an element of $\mathcal{N}$ and an element of $\mathcal{H}$.\\
We now use our extra generator $FUM$ to identify the Klein group $V_4$ inside the group $<G_1,G_2,FUM>$. The group element $(FUM)^3$ is an element of order $2$ as the order of $FUM$ is $6$. So is any conjugate of $(FUM)^3$. When conjugating by $G_2$, we obtain again $(FUM)^3$. But when conjugating by $G_1$, we get a new element of order $2$. Further, the matrix $G_1(FUM)^3\,G_1^{-1}$ commutes to the matrix $(FUM)^3$. Thus, we have
$$\mathcal{V}\,=\,\,<I_3,\; (FUM)^3,\;G_1(FUM)^3\,G_1^{-1},\;(FUM)^3\,G_1\,(FUM)^3\,G_1^{-1}>\;\simeq \;V_4$$
with $V_4 =\mathbb{Z}/2\mathbb{Z}\times\mathbb{Z}/2\mathbb{Z}$.\\
Furthermore, we have the following conjugation relations
\begin{eqnarray}
G_2G_1(FUM)^3\,G_1^{-1}G_2^{-1}&=&(FUM)^3\,G_1(FUM)^3G_1^{-1}\\
G_2^{-1}G_1(FUM)^3\,G_1^{-1}G_2&=&(FUM)^3\,G_1(FUM)^3G_1^{-1}\\
FUM\,G_1(FUM)^3\,G_1^{-1}(FUM)^{-1}&=&G_1(FUM)^3\,G_1^{-1}\\
(FUM)^{-1}\,G_1(FUM)^3\,G_1^{-1}FUM&=&G_1(FUM)^3\,G_1^{-1}\\
G_2(FUM)^3G_1(FUM)^3G_1^{-1}G_2^{-1}&=&G_1(FUM)^3\,G_1^{-1}\\
G_2^{-1}(FUM)^3G_1(FUM)^3G_1^{-1}G_2&=&G_1(FUM)^3\,G_1^{-1}\\
(FUM)^4G_1(FUM)^3G_1^{-1}(FUM)^{-1}&=&(FUM)^3\,G_1(FUM)^3G_1^{-1}\\
(FUM)^2G_1(FUM)^3G_1^{-1}FUM&=&(FUM)^3\,G_1(FUM)^3G_1^{-1}
\end{eqnarray}
%\begin{center} *** These relations are to be understood better. For now they still appear truly magic. The writing may be different here. A conjugation by $(FUM)^3$ swaps the first and last rows and the first and last columns. *** \end{center}

\noindent which show that $\mathcal{V}$ is a normal subgroup of $<G_1,G_2,FUM>$. Some conjugation relations were omitted, which follow from the fact that $(FUM)^3\,G_1(FUM)^3$ is a diagonal matrix. \\
Because the elements of $<A,B>$ have an odd order, while those of $\mathcal{V}$ all have order $2$, we have
$$<A,B>\,\cap\,\mathcal{V}=\lbrace I_3\rbrace$$
Also, we could check with the help of Mathematica that $A$ commutes to both $(FUM)^3$ and $G_1(FUM)^3G_1^{-1}$ and so does the matrix $B$. Then
we have an isomorphism of groups
$$<A,B>\times\mathcal{V}\;\;\simeq\;\; <A,B>\,.\,\mathcal{V},$$
which allows to identify $\mathbb{Z}/18\mathbb{Z}\times\mathbb{Z}/6\mathbb{Z}$ with a subgroup of $Fr(162\times 4)$. \\
Notice further pleasantly that
\begin{eqnarray}
(FUM)\,A\,(FUM)^{-1}&=&A\\
%FUM^{-1}\,A(FUM)&=&A\\
(FUM)\,B(FUM)^{-1}&=&B
%FUM^{-1}B\,(FUM)&=&B
\end{eqnarray}
\noindent thus showing that $<A,B>\lhd\,<G_1,G_2,FUM>$. Then, $$<A,B>.\mathcal{V}\;\lhd\;<G_1,G_2,FUM>$$
We now show that $<A,B>.\mathcal{V}$ intersects $\mathcal{H}$ trivially. %For that, it suffices to check that $H_0^{-1}\,A^i\,B^j\not\in\mathcal{V}$ for all $i$ and $j$ with
%$0\leq i\leq 8$ and $0\leq j\leq 2$. The $27$ products were computed using Mathematica, thus providing the result.
First, if $H_k$ has order $3$, that is $k=3,4$, then since all the elements of $\mathcal{V}$ have order $2$ and commute to $<A,B>$, the matrix $H_k$ would have to be one of the $8$ elements of order $3$ of $\mathcal{N}=<A,B>$. But we know from \cite{BL} that this is impossible as the intersection $\mathcal{N}\cap\mathcal{H}$ is trivial. Next, if $k=0,1,2$, we claim that it suffices to check that
$$\begin{array}{l}(FUM)^3\,H_k\not\in\mathcal{N}\\
G_1(FUM)^3\,G_1^{-1}\,H_k\not\in\mathcal{N}\end{array}$$
when $k=1,2$. When computing the eigenvalues of $(FUM)^3\,H_1$ (resp $(FUM)^3\,H_2$) with the help of Mathematica, one finds out that these are $i$, $-i$ and $1$ (resp $-1$ and $1$) and so the matrix has order $4$ (resp $2$) and cannot belong to $\mathcal{N}$. And when computing the eigenvalues of $G_1\,(FUM)^3\,G_1^{-1}\,H_1$ (resp $G_1\,(FUM)^3\,G_1^{-1}\,H_2$), we find these are $-1$ and $1$ (resp $i$, $-i$ and $1$), and so the same conclusion holds.
Gathering $$\left|\begin{array}{l}\begin{array}{ccc}<A,B>\,.\,\mathcal{V}&\lhd&<G_1,G_2,FUM>\\\mathcal{H}&<&<G_1,G_2,FUM>\end{array}\\\\<A,B>\,.\,\mathcal{V}\cap\,\mathcal{H}=\lbrace I_3\rbrace\end{array}\right.,$$
we obtain
$$(<A,B>.\mathcal{V})\rtimes\mathcal{H}\,\simeq\,<A,B>.\mathcal{V}.\mathcal{H}$$
%If we can show that
%$$<G_1,G_2,FUM>\,\subseteq \mathcal{N}\,\mathcal{V}\,\mathcal{H}\;\;,$$
%then it will follow that $$Fr(162\times 4)\;\simeq\;\mathbb{Z}/18\mathbb{Z}\times\mathbb{Z}/6\mathbb{Z}\rtimes\;S_3$$
We already know from
\cite{BL} that $<G_1,G_2>\,\subseteq\,\mathcal{N}\,\mathcal{H}$.\\\\
Notice further that $FUM\in\mathcal{N}\,\mathcal{V}\,\mathcal{H}$. This is easily seen by simply noticing that
$$(FUM)^{-2}=A^3\;\;,$$ and so
\begin{equation}FUM=A^3(FUM)^3\in\mathcal{N}\mathcal{V}\end{equation}
Note in particular that $<G_1,G_2,FUM>\,=\,\Bigg<G_1,G_2,\begin{pmatrix} &&-1\\&-1&\\-1&&\end{pmatrix}\Bigg>$\\
We deduce $$<G_1,G_2,FUM>\;\subseteq\; \mathcal{N}.\,\mathcal{V}.\,\mathcal{H}$$ So,
by the isomorphism above, we get $$Fr(162\times 4)\;\simeq\;\mathbb{Z}/18\mathbb{Z}\times\mathbb{Z}/6\mathbb{Z}\rtimes\;S_3$$

%\newtheorem{Claim}{Claim}
%\begin{Claim}
%It remains to show that $FUM\in\mathcal{N}\,\mathcal{V}\,\mathcal{H}$.
%\end{Claim}
%\noindent The fact from the claim is straightforward by simply noticing that $$(FUM)^{-2}=A^3\;\;,$$ and so
%$$FUM=A^3(FUM)^3\in\mathcal{N}\mathcal{V}$$
%Note in particular that $<G_1,G_2,FUM>\,=\,\Bigg<G_1,G_2,\begin{pmatrix} &&-1\\&-1&\\-1&&\end{pmatrix}\Bigg>$\\
%\textit{Proof of the Claim.} We also have

\section{Analysis}
In this part, we show that the Freedman group $Fr(162\times 4)$ is isomorphic to the group $D(18,1,1;2,1,1)$. Along the way, we prove the fact that
$\mathbb{Z}/18\mathbb{Z}\times\mathbb{Z}/6\mathbb{Z}\rtimes\;S_3$ is one of the $SU(3)$ finite subgroups from the series $(D)$. \\
In \cite{PO2}, the author studies the structure of the series $(D)$. To that aim, he introduces a new set of generators different from the original set of generators which enlightens the structure of the group. In the case of $D(18,1,1;2,1,1)$ however, this new set of generators coincides with the old one as two of the generators from the new set are simply the identity matrices. We recall some material from \cite{GL2} on the series $(D)$, applied to $D(18,1,1;2,1,1)$. According to \cite{PO2}, the group $D(18,1,1;2,1,1)$ is generated by the following matrices. %\\
%So, using the same notations as in \cite{GL2}, we have $D(18,1,1;2,1,1)=<F(18,1,1),F_1,F^{'}_0,E,B>$ with
$$\begin{array}{l}
F(18,1,1)=\begin{pmatrix} e^{\frac{i\pi}{9}}&&\\&e^{\frac{i\pi}{9}}&\\&&e^{\frac{-2i\pi}{9}}\end{pmatrix}\\\\
E=\begin{pmatrix} 0&1&0\\0&0&1\\1&0&0\end{pmatrix},\;\Bt=\begin{pmatrix}
-1&0&0\\0&0&-1\\0&-1&0\end{pmatrix}\end{array},$$
where the original generator $\overset{\sim}{G}(2,1,1)$ from \cite{PO2} has been renamed $\Bt$, like in \cite{GL2}.
The way the matrix $E$ acts by conjugation is by doing the row and column cycle $(132)$. The way the matrix $\Bt$ acts by conjugation is by swapping the second row and third row and swapping the second column and third column. Consider the subgroup $N(18,1,1;2,1,1)$ generated by the three diagonal matrices $F=F(18,1,1)$, its (identical) $E$- or $\Bt$-conjugate $F^{'}=F^{'}(18,1,1)$ and the $E$-conjugate $F^{''}=F^{''}(18,1,1)$ of $F^{'}$, which we write below for clarity.
$$F^{'}(18,1,1)=\begin{pmatrix} e^{\frac{i\pi}{9}}&&\\&e^{\frac{-2i\pi}{9}} &\\&&e^{\frac{i\pi}{9}}\end{pmatrix},\;F^{''}(18,1,1)=\begin{pmatrix}e^{-\frac{2i\pi}{9}} &&\\&e^{\frac{i\pi}{9}} &\\&&e^{\frac{i\pi}{9}} \end{pmatrix}$$
%All the diagonal matrices from $D(18,1,1;2,1,1)$ are obtained this way. \\
%At this stage, we note that the Freedman group $Fr(162\times 4)$ and the group $D(18,1,1)$ from the series $(D)$ are physically distinct. Indeed, notice
Notice
$$E\Bt=(FUM)^3,$$
so that the Freedman group is also generated by the two braid matrices and the product matrix $E\Bt$. Thus, we have
\begin{eqnarray}
Fr(162\times 4)&=&<G_1,G_2,E\Bt>\\
D(18,1,1;2,1,1)&=&<E,\Bt,F>
\end{eqnarray}
%\begin{eqnarray} Fr(162\times 4)&=&<G_1,G_2,EB>\\
%D(18,1,1;2,1,1)&=&<E,B,F>\end{eqnarray}
%If the braid matrix $G_1$ belonged to the group from the series $D$ with generators originally given back in $1916$ in \cite{MB}, it would have to be one of the diagonal matrices of $<F,F^{'},F^{''}>$.
An algorithm to produce the group $N(18,1,1;2,1,1)$ is provided below.
%Recall $G_1$ is the diagonal matrix with entries $e^{\frac{7i\pi}{9}},\,e^{\frac{13i\pi}{9}},\,e^{\frac{16i\pi}{9}}$.
We start from the triple $(k_1,k_2,k_3)=(1,1,-2)$ corresponding to the $k_i$-th powers of $e^{\frac{i\pi}{9}}$ on the diagonal of $F$. We then build a graph by computing all the possible triples obtained by multiple operations consisting of adding $1$ to two elements of the triple and $-2$ to the remaining one, or adding $-1$ to two elements of the triple and $2$ to the remaining one. It is somewhat like building three parallel interdependent pseudo-random walks in $\mathbb{Z}$, two of which start in $1$ and the last one starts in $-2$ and moving along them in the way described above. In terms of graph of triples, here is what we get.
\begin{center}
\epsfig{file=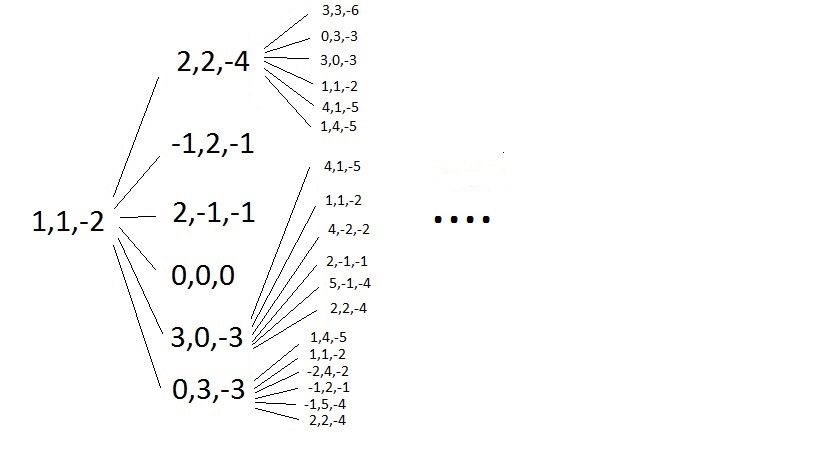, height=10cm}
\end{center}
Of course the order of the group could be large, hence the incomplete graph above could be of large size. But this issue is solved by simply noticing that  $(F^{'})^{-1}$ and $(F^{''})^{-1}$ are both in the second column of the graph above, in respective position $2$ and $3$. This can also be read as
\begin{eqnarray}
F^{'}&=&(F\,F^{''})^{-1}\\
F^{''}&=&(F\,F^{'})^{-1}
\end{eqnarray}
Further, still in the second column of the graph above, two of the triples each form a diagonal matrix of order $6$, namely the triples $(3,0,-3)$ and $(0,3,-3)$.
Consider one of these two matrices, say $F\,(F^{''})^{-1}$ without any loss of generality. From Equality $(16)$, we derive
$$<F,\,F(F^{''})^{-1}>\;=\;<F,F^{'},F^{''}>,$$ as we can write \begin{equation}F^{'}=F^{-2}\,F(F^{''})^{-1}\end{equation}
Furthermore,
because the first two entries in the diagonal matrix corresponding to the triple $(1,1,-2)$ are the same, we obviously have
$$<F>\;\cap\;<F\,(F^{''})^{-1}>\;=\lbrace I_3\rbrace$$
All the diagonal matrices commute, so we derive that
$$N(18,1,1;2,1,1)\,=\;<F,F^{'},F^{''}>\;\simeq\;\mathbb{Z}/18\mathbb{Z}\times\mathbb{Z}/6\mathbb{Z}$$
%and it does not seem reasonable to proceed by hand. Thus, following the algorithm just described, we wrote a program in Mathematica which computes all the distinct such triples, lists them and counts them. We are then able to ask the program whether the triple $(7,13,16)$ corresponding to $G_1$ belongs or not to the list. We found out that the order
Next, as a matter of fact, we have
$$\begin{array}{l}
N(18,1,1;2,1,1)\;\lhd\;D(18,1,1;2,1,1)\\
<E,\Bt>\;\simeq\;S_3\\
<E,\Bt>\;\cap\;N(18,1,1;2,1,1)=\lbrace I_3\rbrace
\end{array}$$
So, $D(18,1,1;2,1,1)=N(18,1,1;2,1,1)\rtimes S_3$, where the action by $S_3$ on the group $N(18,1,1;2,1,1)$ is given by
$$\begin{array}{ccccc}
\Bt\,F\,\Bt^{-1}&=&E\,F\,E^{-1}&=&F^{'}\end{array}$$
$$\begin{array}{cc}\begin{array}{ccc}
\Bt\,F^{'}\,\Bt^{-1}&=&F\\
\Bt\,F^{''}\,\Bt^{-1}&=&F^{''}\end{array}
&\begin{array}{ccc}
E\,F^{'}\,E^{-1}&=&F^{''}\\
E\,F^{''}\,E^{-1}&=&F\end{array}
\end{array}$$
In terms of the generators $F$ and $F(F^{''})^{-1}$, the conjugation relations can easily be read on the "horizontal tree". Indeed, the conjugate $\Bt\,F\,(F^{''})^{-1}\,\Bt^{-1}$ is the triple $(3,-3,0)$, that is a descendant of $(2,-1,-1)$ in the incomplete third column of the tree,
and the conjugate $E\,F\,(F^{''})^{-1}\,E^{-1}$ is the triple $(0,-3,3)$, whose inverse appears at the very bottom of the second column in the tree. \\We then deduce
$$\begin{array}{ccccccc}
\Bt\,F\,(F^{''})^{-1}\,\Bt^{-1}&=&F(F^{'})^2&\overset{1}{=}&F\,F^{-4}[F\,(F^{''})^{-1}]^2&=&F^{-3}[F(F^{''})^{-1}]^2\\
E\,F\,(F^{''})^{-1}\,E^{-1}&=&(F(F^{'})^{-1})^{-1}&=&F^{-1}F^{'}&\overset{2}{=}&F^{-3}\,F(F^{''})^{-1}
\end{array}$$
Equalities $1$ and $2$ in the series of equalities above hold by Eq. $(18)$. \\
In summary, the group $D(18,1,1;2,1,1)$ is
$$\begin{array}{l}\bigg<\begin{pmatrix} e^{\frac{i\pi}{9}}&&\\&e^{\frac{i\pi}{9}}&\\&&e^{-2\frac{i\pi}{9}}\end{pmatrix},\;\begin{pmatrix} e^{\frac{i\pi}{3}}&&\\&1&\\&&e^{-\frac{i\pi}{3}}\end{pmatrix}\;\bigg>\;\;\rtimes\;\;\bigg<\begin{pmatrix} 0&1&0\\0&0&1\\1&0&0\end{pmatrix},\;
\begin{pmatrix} -1&0&0\\0&0&-1\\0&-1&0\end{pmatrix}\bigg>\\\\ \simeq \mathbb{Z}/18\mathbb{Z}\times\mathbb{Z}/6\mathbb{Z}\;\rtimes\;S_3\end{array}$$
\noindent for the action provided above. The set of matrices from the symmetric group
$$S_3=\lbrace id,(23),(12),(13),(132),(123)\rbrace$$ are the following (in the same order)
$$S_3=\bigg\lbrace\begin{array}{l}\begin{pmatrix}1&&\\&1&\\&&1\end{pmatrix},\;\begin{pmatrix}-1&0&0\\0&0&-1\\0&-1&0\end{pmatrix},\;\begin{pmatrix}
0&-1&0\\-1&0&0\\0&0&-1\end{pmatrix},\\\begin{pmatrix}&&-1\\&-1&\\-1&&\end{pmatrix},\;\begin{pmatrix}0&1&0\\0&0&1\\1&0&0\end{pmatrix},\;\begin{pmatrix}
0&0&1\\1&0&0\\0&1&0\end{pmatrix}\end{array}\bigg\rbrace,$$

\noindent that is (still in the same order)
$$S_3=\lbrace I_3,\,\Bt,\,\Bt E,\,E\Bt,\,E,\,\Bt E\Bt\rbrace$$
Let us now go back to our first group, the Freedman group $Fr(162\times 4)$. We must provide a generator $C_{18}$ for the cyclic group $\mathbb{Z}_{18}$ which will be given by the product of the element $A$ of order $9$ with an element of order $2$ from the Klein group, say $(FUM)^3$, the cube of the Freedman generator.
$$C_{18}=\,A\,(FUM)^3=\;e^{\frac{5i\pi}{9}}\;\begin{pmatrix} -\frac{1}{2}(1+e^{\frac{i\pi}{3}})&0&-\frac{1}{2}(-1+e^{\frac{i\pi}{3}})\\
0&1&0\\
-\frac{1}{2}(-1+e^{\frac{i\pi}{3}})&0&-\frac{1}{2}(1+e^{\frac{i\pi}{3}})\end{pmatrix}$$
\noindent And likewise do so for the cyclic group $\mathbb{Z}_6$ which will be generated by $B$ times a different element of order two from the Klein group, say
$G_1(FUM)^3G_1^{-1}$.
It yields $$C_6=\,B\,G_1\,(FUM)^3\,G_1^{-1}=\begin{pmatrix}\frac{1}{2}(-1+e^{\frac{2i\pi}{3}})&0&\frac{1}{2}(1+e^{\frac{2i\pi}{3}})\\
0&e^{\frac{i\pi}{3}}&0\\\frac{1}{2}(1+e^{\frac{2i\pi}{3}})&0&\frac{1}{2}(-1+e^{\frac{2i\pi}{3}}) \end{pmatrix}$$
As far as the symmetric group, recall the elements of order two are the three matrices $H_0$, $H_1$ and $H_2$ introduced earlier and those of order three are the two matrices $H_3$ and $H_4$, also defined earlier. Like for $D(18,1,1;2,1,1)$, we can write
$$S_3=\lbrace I_3,\,H_1,\,H_1H_3,\,H_3H_1,\,H_3,\,H_1H_3H_1\rbrace$$ Then, a presentation of the Freedman group is as stated in Theorem $1$. Here is how we obtain the conjugation relations.
First, we notice the equality
\begin{equation}H_3\,C_{18}H_3^{-1}=C_{18}^{10}C_6^4\begin{pmatrix}&&-1\\&-1&\\-1&&\end{pmatrix}\begin{pmatrix}-1&&\\&1&\\&&-1\end{pmatrix}\end{equation}
But \begin{equation} C_6^3=\begin{pmatrix}&&1\\&-1&\\1&&\end{pmatrix}\end{equation}
%from which we derive another equality \begin{equation}H_3,C_{18}H_3^{-1}=C_{18}^7\,C_6^4\,\begin{pmatrix}&&-e^{\frac{2i\pi}{3}}\\&-e^{\frac{2i\pi}{3}}&\\-e^{\frac{2i\pi}{3}}&&\end{pmatrix}=C_{18}^7\,C_6^4\,FUM\end{equation}
%But recall from earlier in $(13)$ that $FUM=A^3\,(FUM)^3$, so that
%\begin{equation}FUM=A^2\,C_{18}\end{equation}
%Further, by Eq. $(11)$, we know that $FUM$ commutes to $A$. And $FUM$ has order $6$. Therefore,
%\begin{equation}
%C_{18}^2=A^2
%\end{equation}
%Combining Eq. $(19)$, $(20)$ and $(21)$ yields
Hence we get \begin{equation}H_3C_{18}H_3^{-1}=C_{18}^{10}\,C_6^7\,,\end{equation}
\\like in the statement of Theorem $1$. \\\\Note for future references%along the way that Eq. $(20)$ and $(21)$ imply in particular
\begin{equation}
C_{18}^2=A^2
\end{equation}
and \begin{equation}
C_{18}^3=FUM
\end{equation}
We notice that
\begin{equation}
H_1\,C_6\,H_1^{-1}=C_6^2(FUM)^2\begin{pmatrix} &&1\\&-1&\\1&&\end{pmatrix}
\end{equation}
%We also notice that
%\begin{equation}
%C_6^3\,FUM=\begin{pmatrix} e^{-\frac{i\pi}{3}}&&\\&e^{\frac{2i\pi}{3}}&\\&&e^{-\frac{i\pi}{3}}\end{pmatrix}
%\end{equation}
%Since $$(FUM)^2=\begin{pmatrix} -e^{\frac{i\pi}{3}}&&\\&-e^{\frac{i\pi}{3}}&\\&&-e^{\frac{i\pi}{3}}\end{pmatrix}\,,$$
%we get
%$$(FUM)^2\,C_6^3\,FUM=\begin{pmatrix} -1&&\\&1&\\&&-1\end{pmatrix}=(FUM)^3\,\begin{pmatrix}&&1\\&-1&\\1&&\end{pmatrix}$$
%Then,
%\begin{equation}
%\begin{pmatrix}&&1\\&-1&\\1&&\end{pmatrix}=(FUM)^{-1}\,C_6^3\,FUM
%\end{equation}
%By Eq. $(23)$, $(24)$ and $(26)$, it now comes
Hence \begin{equation}
H_1\,C_6\,H_1^{-1}=C_6^5\,C_{18}^6
\end{equation}
Finally, we have
\begin{equation}
(H_3\,C_6\,H_3^{-1})(FUM\,C_{18}^6\,C_6^5)=\begin{pmatrix}e^{-\frac{2i\pi}{3}}&&\\&e^{-\frac{2i\pi}{3}}&\\&&e^{-\frac{2i\pi}{3}}\end{pmatrix}
\end{equation}
\begin{equation}
(H_1C_{18}H_1^{-1})(FUM\,C_{18}^4\,C_6^5)=\begin{pmatrix}e^{\frac{4i\pi}{9}}&&\\&-e^{\frac{i\pi}{9}}&\\&&e^{\frac{4i\pi}{9}}\end{pmatrix}
\end{equation}
In order to find Eq. $(26)$ and $(27)$, we used the fact that
\begin{equation}\begin{pmatrix} a&0&b\\0&c&0\\b&0&a\end{pmatrix}\begin{pmatrix}a &0&-b\\0&c&0\\-b&0&a\end{pmatrix}=\begin{pmatrix} a^2-b^2&&\\&c^2&\\&&a^2-b^2\end{pmatrix}\end{equation}
after we noticed that the respective right hand side factors in the products above differed from their respective left hand side factors
by two minus signs on the two extremities of the anti-diagonal.
The matrix in Eq. $(26)$ is $(FUM)^2$, hence we get after using Eq. $(23)$ and simplifying
\begin{equation}
H_3\,C_6\,H_3^{-1}=C_6^{-5}\,C_{18}^{-3}
\end{equation}
After investigation, the matrix in Eq. $(27)$ is \begin{equation}C_6^2\,C_{18}^8\end{equation}
We thus obtain after simplification
\begin{equation}
H_1\,C_{18}\,H_1^{-1}=C_6^3\,C_{18}
\end{equation}

\noindent Eqs. $(21)$, $(25)$, $(29)$ and $(31)$ provide the conjugation relations of Theorem $1$.
After a quick glance at the presentation of Theorem $1$, this is not the same presentation as the one given earlier for $D(18,1,1;2,1,1)$. Then, we must still investigate whether there could be an isomorphism between the two semi-direct products
$$\mathbb{Z}/18\mathbb{Z}\times\mathbb{Z}/6\mathbb{Z}\underset{\varphi_{_D}}{\rtimes}S_3\;\;\;\text{and}\;\;\;\mathbb{Z}/18\mathbb{Z}\times
\mathbb{Z}/6\mathbb{Z}\underset{\varphi_{_F}}{\rtimes}S_3$$ where
$$\begin{array}{cccc}
\varphi_{_D}:&S_3&\longrightarrow&Aut(\mathbb{Z}_{18}\times\mathbb{Z}_6)\\
\varphi_{_F}:&S_3&\longrightarrow&Aut(\mathbb{Z}_{18}\times\mathbb{Z}_6)
\end{array}$$
are the two homomorphisms provided in this paper for the respective groups $D(18,1,1;2,1,1)$ and $Fr(162\times 4)$.
%If such an isomorphism were to exist, it would have to map the unique $3$-Sylow subgroup of $D(18,1,1;2,1,1)$ onto the unique $3$-Sylow subgroup of $Fr(162\times 4)$.
First, we clearly state that the groups $Fr(162)$ and $D(9,1,1;2,1,1)$ are isomorphic groups like claimed in
\cite{BL} and mentioned several times along this paper. We provide an explicit isomorphism between both groups.
\begin{Theorem}
The groups $Fr(162)$ and $D(9,1,1;2,1,1)$ are isomorphic groups.
\end{Theorem}
First we establish the following result.
\newtheorem{Lemma}{Lemma}
\begin{Theorem}
The unique $3$-Sylow subgroup of $Fr(162)$ is
\begin{equation}\mathcal{S}_3(F)=\mathcal{N}\,\bigsqcup\,\mathcal{N}H_3\,\bigsqcup\,\mathcal{N}H_3^2\end{equation}\end{Theorem}

The fact that the unions are disjoint comes from the following facts below. A proof of Eq. $(32)$ is then given.
\newtheorem{Definition}{Definition}
\begin{Definition}
Call "cross matrix" a $3$ by $3$ matrix having zeroes in positions $(i,j)$ where exactly one of $i$ or $j$ is even.
\end{Definition}
\begin{Lemma}
The matrices $A$ and $B$ are cross matrices and the set of special unitary cross matrices forms a subgroup of $SU(3)$.
\end{Lemma}
The proof of the second point of Lemma $1$ is straightforward since the inverse of a unitary matrix is its conjugate transpose.
\newtheorem{Corollary}{Corollary}
\begin{Corollary}
All the matrices $A^i\,B^j$ with $0\leq i\leq 8$ and $0\leq j\leq 3$ are cross matrices.
\end{Corollary}
\noindent Since $H_3$ and $H_3^2$ are not cross matrices, by the corollary, the unions above are indeed disjoint.
This provides $27+27+27=3^4$ elements. We now finish proving that Eq. $(32)$ holds.
\begin{Lemma}
Let $\mathcal{S}_3$ denote the unique $3$-Sylow subgroup of $Fr(162)$. If $x$ does not belong to $\mathcal{S}_3$, then $x^2$ belongs to $\mathcal{S}_3$.
\end{Lemma}
\textsc{Proof.} Since $162=3^4\times 2$, the number $n_3$ of $3$-Sylow subgroups of $Fr(162)$ satisfies to
$$\begin{array}{l}
n_3|2\\
n_3\equiv 1\;(mod 3)
\end{array}$$
Both conditions imply that $n_3=1$. Hence $Fr(162)$ contains a unique $3$-Sylow subgroup, say $\mathcal{S}_3$, which is normal in $Fr(162)$.
It now suffices to consider the quotient group $Fr(162)/\mathcal{S}_3$. This is a group of order $2$. If an element $x$ of $Fr(162)$ does not belong to the $3$-Sylow, then its image in the quotient group has order $2$. That is $x^2$ belongs to the $3$-Sylow.
\begin{Corollary}
Both $H_3$ and $H_3^2$ belong to the $3$-Sylow.
\end{Corollary}
\textsc{Proof.} By contradiction, if $H_3\not\in\mathcal{S}_3$, then by the lemma above, we must have $H_3^2\in\mathcal{S}_3$. Then both elements are in fact in the $3$-Sylow since $H_3^2$ has order $3$.\\

\textit{End of the proof of Theorem $3$.} If we can show that $\mathcal{S}_3$ contains $\mathcal{N}$, we are done. But the fact that $B$ belongs to $\mathcal{S}_3$ follows from Lemma $2$ by the same argument as before. And the fact that $A$ belongs to $\mathcal{S}_3$ is also a consequence of Lemma $2$. Indeed, by contradiction, if $A$ does not belong to $\mathcal{S}_3$, then $A^2$ belongs to $\mathcal{S}_3$. But because the integer $2$ is prime to $9$, the element $A^2$ also generates the cyclic group $<A>$ of order $9$. \\
The fact that $\mathcal{N}$ and $<H_3>=\lbrace I_3,H_3,H_3A^2\rbrace$ are both contained in $\mathcal{S}_3$ and the count from above suffice to imply that the disjoint union of the right hand side of Eq. $(32)$ is a group which is the unique $3$-Sylow of $Fr(162)$.\\\\
We proceed by now determining the unique $3$-Sylow $\mathcal{S}_3(D)$ of $D(9,1,1;2,1,1)$. The group $D(9,1,1;2,1,1)$ is the group generated by the three matrices $E$, $\Bt$ and $F^2$, where we used the same notations as before. The abelian group of $D(18,1,1;2,1,1)$ composed of diagonal matrices contains a subgroup, say $\mathcal{F}$ with $$\mathcal{F}\,=\,<F^2>\times <[F(F^{''})^{-1}]^2>$$
This group is isomorphic to $\mathbb{Z}/9\mathbb{Z}\times\mathbb{Z}/3\mathbb{Z}$. It is also a subgroup of $D(9,1,1;2,1,1)$.
\begin{Lemma}
The unique $3$-Sylow subgroup $\mathcal{S}_3(D)$ of $D(9,1,1;2,1,1)$ must contain all the elements of odd order.
\end{Lemma}
Indeed, suppose $x$ is an element of odd order of $D(9,1,1;2,1,1)$ which is not in $\mathcal{S}_3(D)$. Then, by considering the quotient, we must have $x^2$ belongs to $\mathcal{S}_3(D)$. But since $x$ has odd order, the element $x^2$ is also a generator of the cyclic group $<x>$. This represents a contradiction.
As a corollary to the lemma, the direct product $\mathcal{F}$ is contained in
the $3$-Sylow $\mathcal{S}_3(D)$. And the two matrices $E$ and $E^2$ of order $3$ also belong to $\mathcal{S}_3(D)$. Moreover, these two matrices are not diagonal. Thus, we conclude
\begin{Theorem} The unique $3$-Sylow subgroup of $D(9,1,1;2,1,1)$ is
\begin{equation}\mathcal{S}_3(D)=\mathcal{F}\,\bigsqcup\,\mathcal{F}\,E\,\bigsqcup\,\mathcal{F}\,E^2\end{equation}
Moreover, this group is $C(9,1,1)$.
\end{Theorem}
From the previous study, we derive some facts about the structure of the bigger groups $D(18,1,1;2,1,1)$ and $Fr(162\times 4)$.
\begin{Theorem}
Both $D(18,1,1;2,1,1)$ and $Fr(162\times 4)$ contain exactly four $3$-Sylows.
\end{Theorem}
\textsc{Proof.} For each group, the number of $3$-Sylows is either $1$ or $4$. Suppose for a contradiction that this number is one. The unique $3$-Sylow subgroups would have to be $\mathcal{S}_3(D)$ and $\mathcal{S}_3(F)$ respectively. Also, they would have to contain all the elements of odd order of their respective groups. Indeed if an element $x$ does not belong to the unique $3$-Sylows, then by arguments similar as before, $x^2$ or $x^4$ or $x^8$ belongs to the $3$-Sylow. When $x$ has an odd order, any of these elements is also a generator for the cyclic group $<x>$, hence a contradiction.
As far as $D(9,1,1;2,1,1)$, simply notice the product matrix $F^{'}E$ has order $3$ and does not belong to $\mathcal{S}_3(D)$ since
$$F^{'}=F^{-2}\,F(F^{''})^{-1}\not\in\mathcal{F}$$
And in the case of $Fr(162\times 4)$, the element $(FUM)\,H_3$ has order $3$ and does not belong to $\mathcal{S}_3(F)$, as neither of $(FUM)\,H_3$ nor $(FUM)\,H_3^{-1}$ is a cross matrix.\\\\

In light of Theorem $3$ and $4$, it is now easy to conclude that $D(9,1,1;2,1,1)$ and $Fr(162)$ are isomorphic, like stated in Theorem $2$. An isomorphism between $Fr(162)$ and $D(9,1,1;2,1,1)$ must map $\mathcal{S}_3(F)$ onto $\mathcal{S}_3(D)$.
\begin{Lemma}\hfill\\\\
$(i)$ The only elements of order $3$ which belong to the center of the unique $3$-Sylow $\mathcal{S}_3(F)$ of $Fr(162)$ are $A^3$ and $A^6$.\\\\
$(ii)$ The only elements of order $3$ which belong to the center of the unique $3$-Sylow $\mathcal{S}_3(D)$ of $D(9,1,1;2,1,1)$ are $F^6$ and $F^{12}$.
\end{Lemma}
\textsc{Proof.} $(i)$ The matrices $H_3$ and $B$ don't commute. Likewise, the matrices $H_3^2$ and $B$ don't commute. So an element of order $3$ belonging to $Z(\mathcal{S}_3(F))$ would have to belong to $\mathcal{N}$. The elements of order $3$ of $\mathcal{N}$ are
$$A^3,A^6,B,B^2,A^3\,B,A^3\,B^2,A^6\,B,A^6\,B^2$$
Among these, only $A^3$ and $A^6$ do commute to $H_3$.\\
$(ii)$ The matrices $E$ (resp $E^2$) and $F^2$ don't commute. Hence an element of order $3$ belonging to $Z(\mathcal{S}_3(D))$ must belong to $\mathcal{F}$.
The fact that $F^6$ commutes to $E$ while $[F(F^{''})^{-1}]^2$ does not commute to $E$ implies the result. \\

A corollary to this Lemma is that if $f$ is an isomorphism between $Fr(162)$ and $D(9,1,1;2,1,1)$, then
$$\left\lbrace\begin{array}{l}
f(A^3)=F^6\;\text{and}\;f(A^6)=F^{12}\\
\text{or}\\
f(A^3)=F^{12}\;\text{and}\;f(A^6)=F^6
\end{array}\right.$$
It gives some ingredients for the following Theorem.
\begin{Theorem}
The map
$$\begin{array}{ccc}
Fr(162)&&D(9,1,1;2,1,1)\\
&\overset{\sim}{\longrightarrow}&\\
<A,B>\rtimes <H_1,H_3> & & <F^2,[F(F^{''})^{-1}]^2>\rtimes <\overset{\sim}{B},E>\\
&&\\
&&\\
A&\longmapsto&F^2\\
B&\longmapsto&F^{12}[F(F^{''})^{-1}]^2\\
H_1&\longmapsto&\overset{\sim}{B}\,E\\
H_3&\longmapsto&E
\end{array}$$
defines an isomorphism of groups between $Fr(162)$ and $D(9,1,1;2,1,1)$
\end{Theorem}
\textsc{Proof.} It suffices to show that the following relations of the presentation of $Fr(162)$
\begin{eqnarray}
H_1\,A\,H_1^{-1}&=&A\\
H_1\,B\,H_1^{-1}&=&A^6\,B^2\\
H_3\,A\,H_3^{-1}&=&A\,B\\
H_3\,B\,H_3^{-1}&=&A^6\,B\,
\end{eqnarray}
derived from the end of \cite{BL} with $H_1=G_1\,G_2\,G_1=T_1$, $H_0=T_3$ and $H_3=T_1T_3$ are satisfied on the images.
As a matter of fact, we have
\begin{equation}
\overset{\sim}{B}\,E\,F^2\,E^{-1}\,\overset{\sim}{B}=F^2
\end{equation}
which motivated our choice for the image of $H_1$.
Next we must have
\begin{equation}
f(H_3)F^2f(H_3)^{-1}=F^2\,f(B)\;,
\end{equation}
from which we derive
\begin{equation}
f(B)=F^{-2}\,f(H_3)F^2\,f(H_3)^{-1}
\end{equation}
Choose to map $H_3$ to $E$ by $f$. Then, it follows
\begin{equation}
f(B)=F^{-2}\,EF^2E^{-1}=\begin{pmatrix}1&&\\&-e^{\frac{i\pi}{3}}&\\&&e^{\frac{2i\pi}{3}}\end{pmatrix}=F^{12}\,[F(F^{''})^{-1}]^2
\end{equation}
We checked with Mathematica that
\begin{eqnarray}
\overset{\sim}{B}EF^{-2}EF^2E^{-1}E^{-1}\overset{\sim}{B}&=&F^{12}F^{-2}EF^2E^{-1}F^{-2}EF^2E^{-1}\\
EF^{-2}EF^2E^{-1}E^{-1}&=&F^{12}F^{-2}EF^2E^{-1}
\end{eqnarray}
Hence Eqs. $(35)$ and $(37)$ are satisfied on the images, which ends the proof of Theorem $6$.\\\\

We now continue our study of the group extensions. It will be useful to explicit the four $3$-Sylow subgroups. A first result is as follows.
\begin{Lemma}
The four $3$-Sylows of $Fr(162\times 4)$ are the following.\\
$$\begin{array}{cc}\begin{array}{cc}(i)& \mathcal{S}_3(F)\\
(ii)& V_2\;\mathcal{S}_3(F)\,V_2^{-1}\\
(iii)& V_3\;\mathcal{S}_3(F)\,V_3^{-1}\\
(iv)& V_4\;\mathcal{S}_3(F)\,V_4^{-1}
\end{array}&,\end{array}$$
where the elements of the Klein group $\mathcal{V}$ have been renamed $V_2$, $V_3$ and $V_4$ for convenience. For instance, $V_2=(FUM)^3$.
\end{Lemma}
We have the following corollary.
\begin{Corollary}
Denote by $3SylF$ the subgroup of $Fr(162\times 4)$ generated by the four $3$-Sylows. \\We have
$$3SylF\;=\;<\mathcal{V},\mathcal{S}_3(F)>\;\simeq \mathcal{V}\rtimes\mathcal{S}_3(F)$$
Consequently, $Fr(162\times 4)$ is not generated by its $3$-Sylow subgroups. In other words, $3SylF$ is a proper normal subgroup of $Fr(162\times 4)$.
We have the series
$$\mathcal{N}.\mathcal{V}\;\lhd\;3SylF\;\lhd\;Fr(162\times 4),$$
where the quotients are simple.
\end{Corollary}
%\begin{Lemma}
%Let $G$ be a finite group and $p$ a prime number. \\\\
%$(i)$ The subgroup generated by all the $p$-Sylows of $G$ is normal in $G$. \\
%$(ii)$ The group $G$ is generated by its $p$-Sylows if and only if the order of any non-trivial quotient of $G$ is divisible by $p$.
%\end{Lemma}
%\begin{Theorem}
%The Freedman group $Fr(162\times 4)$ is generated by its $3$-Sylow subgroups.
%\end{Theorem}
%\begin{Theorem}
%We have the following two subnormal series.

%\begin{center}$$\begin{array}{ccccccc} \mathcal{S}_3(D)&\lhd&D(9,1,1;2,1,1)&\lhd& <3\text{ Sylows of}\;D(18,1,1;2,1,1)>&\lhd& D\\
%&&&&&&\\
%|\Large{S}&&|\Large{S}&&&&\\
%&&&&&&\\
%\mathcal{S}_3(F)&\lhd&Fr(162)&\lhd& <3\text{ Sylows of}\;Fr(162\times 4)>&\lhd&Fr
%\end{array}$$\end{center}
%\end{Theorem}
%\noindent In the Theorem, $D$ denotes $D(18,1,1;2,1,1)$ and $Fr$ denotes the Freedman group $Fr(162\times 4)$. The brackets denote the respective subgroups generated by all the $3$-Sylows inside $D$ (resp $Fr$).\\\\
\noindent \textsc{Proof.} %First, we show that $D(9,1,1;2,1,1)$ is a subgroup of the group $3SylF$ generated by the four $3$-Sylows of $D(18,1,1;2,1,1)$.
By Theorem $3$, we know that
$$\mathcal{S}_3(F)=\mathcal{N}\bigsqcup\mathcal{N}H_3\bigsqcup\mathcal{N}H_4$$ is one of the four $3$-Sylows of $D$. Another distinct one is
$$(FUM)\mathcal{S}_3(F)(FUM)^{-1}=\mathcal{N}\bigsqcup\mathcal{N}[(FUM)\,H_3(FUM)^{-1}]\bigsqcup\mathcal{N}[(FUM)\,H_4(FUM)^{-1}]$$
since both $A$ and $B$ are self-conjugate under $FUM$. %To show the desired inclusion, it suffices to show that $H_0$ belongs to $3SylF$. Recall
%$$H_0=\begin{pmatrix} -1&&\\ &-1&\\&&1\end{pmatrix}$$
We have
$$H_4(FUM)H_3(FUM)^{-1}=\begin{pmatrix} 0&0&1\\0&-1&0\\1&0&0\end{pmatrix}\in 3SylF$$
and
$$(FUM)H_3(FUM)^{-1}H_4=\begin{pmatrix} -1&0&0\\0&1&0\\0&0&-1\end{pmatrix}\in 3SylF$$
The product of these two matrices is $$\begin{pmatrix}0&0&-1\\0&-1&0\\-1&0&0\end{pmatrix}$$
So, $(FUM)^3$ belongs to $3SylF$ (and so does $FUM$ by Eq. $(13)$). Then, by the way it is defined, the whole Klein group $\mathcal{V}$ is contained in $3SylF$.
This implies
\begin{equation}
<\mathcal{V},\mathcal{S}_3(F)>\;\subseteq\, 3SylF
\end{equation}
The subgroups of Lemma $5$ are all $3$-Sylows. It remains to show that they are all distinct. Notice the products
$$V_jH_3V_j^{-1}H_k$$
with $j\in\lbrace 2,3,4\rbrace$ and $k\in\lbrace 3,4\rbrace$, are either not cross-matrices or have order two or are $(FUM)^3$ or $V_3$. In any case, they do not belong to $\mathcal{N}$. This shows the result. Lemma $5$ implies that
\begin{equation}
3SylF\;\subseteq\;<\mathcal{V},\mathcal{S}_3(F)>
\end{equation}
Gathering the two inclusions $(44)$ and $(45)$, we obtain
\begin{equation}
3SylF\;=\;<\mathcal{V},\mathcal{S}_3(F)>\; ,
\end{equation}
as stated in Corollary $3$. Recall that $\mathcal{V}$ is a normal subgroup of $Fr(162\times 4)$. From $(46)$, we then derive
\begin{equation}
3SylF/\mathcal{V}\;\simeq \;\mathcal{S}_3(F)
\end{equation}
It follows that \begin{equation}|3SylF|=\;4.3^4\end{equation} by Lagrange Theorem. In particular, $3SylF$ is a proper subgroup of $Fr(162\times 4)$. And the fact that it is a normal subgroup is a straightforward consequence of the fact that all the $3$-Sylows are conjugate or follows from the fact that $3SylF$ has index two in $Fr(162\times 4)$. Finally, since
$$\begin{array}{l}
\mathcal{V}\cap\mathcal{S}_3(F)=\lbrace I_3\rbrace\\
|<\mathcal{V},\mathcal{S}_3(F)>|=|\mathcal{V}||\mathcal{S}_3(F)|\;,\\
\mathcal{V}\lhd\;Fr(162\times 4)
\end{array}$$
the map
$$\begin{array}{ccc}
\mathcal{V}\rtimes\mathcal{S}_3(F)&\longrightarrow& <\mathcal{V},\mathcal{S}_3(F)>\\
(V_i,NH_k)&\longmapsto&V_iNH_k\\&&
\end{array}$$
with $N$ an element of $\mathcal{N}$ and $k\in\lbrace 3,4\rbrace$ and $V_1=I_3$ and $i\in\lbrace 1,2,3,4\rbrace$,
is an isomorphism of groups which allows to identify the subgroup $3SylF$ of $Fr(162\times 4)$ generated by all the $3$-Sylows of $Fr(162\times 4)$ with the semi-direct product $\mathcal{V}\rtimes\mathcal{S}_3(F)$ with respect to conjugation, where $\mathcal{S}_3(F)$ is the unique $3$-Sylow subgroup of $Fr(162)$ and $\mathcal{V}$ is the Klein group formerly defined.
\begin{Lemma}
The four $3$-Sylows of $D(18,1,1;2,1,1)$ are
$$\begin{array}{cc}(i)& \mathcal{S}_3(D)\\
(ii)& F\;\mathcal{S}_3(D)\,F^{-1}\\
(iii)& EF\;\mathcal{S}_3(D)\;F^{-1}E^{-1} \\
(iv)& E^2F\;\mathcal{S}_3(D)\;F^{-1}E^{-2}
\end{array}$$
\end{Lemma}
First, we show that

$$\begin{array}{l}F\,\mathcal{S}_3(D)\,F^{-1}=\mathcal{F}\;\bigsqcup\; \mathcal{F}(FEF^{-1})\;\bigsqcup\;\mathcal{F}(FE^2F^{-1})\\\\
F\mathcal{S}_3(D)F^{-1}\neq \mathcal{S}_3(D)\\
\end{array}$$

This follows from three simple facts:\\\\
 $1)$ $FEF^{-1}$ is not a diagonal matrix, hence the unions above are disjoint.\\
  $2)$ $FEF^{-1}E^{-1}$ is diagonal but has order $6$, hence does not belong to the group $\mathcal{F}$ of odd order $27$\\
   $3)$ $FEF^{-1}E^{-2}$ is not a diagonal matrix.\\

Next, the groups of $(i)$ and $(iii)$ are distinct since $F^{-1}EFE^{-1}$ has order $6$ and $F^{-1}EFE^{-2}$ is not diagonal. And the groups of $(i)$ and $(iv)$ are distinct since $E^{-1}F^{-1}E^2FE^{-1}$ has order $6$ and $E^{-1}F^{-1}E^2FE^{-2}$ is not diagonal. It remains to show that the two groups in $(ii)$ and $(iii)$ are distinct and so are those of $(ii)$ and $(iv)$. This follows from similar arguments.

\begin{Lemma}Denote by $3SylD$ the subgroup of $D(18,1,1;2,1,1)$ generated by all the $3$-Sylows of $D(18,1,1;2,1,1)$.
The group $3SylD$ contains the Klein group, say $\mathcal{V}_D$,
$$I_3,\;W_2=\begin{pmatrix} -1&&\\&1&\\&&-1\end{pmatrix},\;W_3=\begin{pmatrix}1&&\\&-1&\\&&-1\end{pmatrix},\;W_4=\begin{pmatrix}-1&&\\&-1&\\&&1\end{pmatrix}$$
Moreover, $\mathcal{V}_D$ is normal in $D(18,1,1;2,1,1)$.
\end{Lemma}
\textsc{Proof.} Notice the matrices $FEF^{-1}E^2$ and $E^2FEF^{-1}$ both have order $6$. The cube of the first one (resp second one) is the second (resp third) matrix from the Lemma. The product of their respective cubes is the last matrix of the Lemma.
\begin{Lemma}
The four $3$-Sylows of $D(18,1,1;2,1,1)$ are also given by $\mathcal{S}_3(D)$ and its $W_i$-conjugates.
\end{Lemma}
\textsc{Proof.} It is easy to check by similar arguments as already used before that they are all distinct. Hence the immediate corollary.
\begin{Corollary}
$3SylD\;=\;<\mathcal{V}_D,\mathcal{S}_3(D)>\;\simeq\;\mathcal{V}_D\rtimes\mathcal{S}_3(D)$. Moreover, we have
$$N(18,1,1;2,1,1)=\mathcal{V}_D\;.\;\mathcal{F}\simeq \mathcal{V}_D\times\mathcal{F}$$
and the series
$$\mathcal{V}_D\;.\;\mathcal{F}\;\lhd\;3SylD\;\lhd\;D(18,1,1;2,1,1)\;,$$
where the quotients are simple.
\end{Corollary}
In light of Corollary $3$ and $4$, we are now ready to show the following Theorem.
\begin{Theorem}
The Freedman group is isomorphic to the finite SU(3)-subgroup $D(18,1,1;2,1,1)$ from the series $D$.
\end{Theorem}
\textsc{Proof.} Suppose that there exists an isomorphism of groups, say $g$ between $Fr(162\times 4)$ and $D(18,1,1;2,1,1)$. Then, $g$
maps $$3SylF\;=\;\mathcal{V}\rtimes\mathcal{S}_3(F)$$
onto $$3SylD\;=\;\mathcal{V}_D\rtimes\mathcal{S}_3(D)$$
Further, we claim that such an isomorphism $g$ maps $\mathcal{V}$ onto $\mathcal{V}_D$. Indeed, an element of order $2$ belonging to $\mathcal{V}$ must be mapped to an element of order $2$ of $3SylD$. The elements of $3SylD$ can be uniquely written as a product of an element of $\mathcal{V}_D$ times an element of $\mathcal{S}_3(D)$. The matrices of $\mathcal{V}_D$ are all diagonal matrices, $W_i$'s as we called them earlier. An element of $\mathcal{S}_3(D)$ is a diagonal matrix of $\mathcal{F}$, say $\Delta$, times $I_3$ or $E$ or $E^2$. It is now easy to see that the only elements of order two of $3SylD$ are those belonging to the Klein group $\mathcal{V}_D$. Indeed, we have
\begin{eqnarray}
(W_i\,\Delta\,E)\,(W_i\,\Delta\,E)&=&(W_i\,\Delta)\,(E\,W_i\,\Delta\,E)\\
(W_i\,\Delta\,E^2)\,(W_i\,\Delta\,E^2)&=&(W_i\,\Delta)\,(E^2\,W_i\,\Delta\,E^2)
\end{eqnarray}
The effect of multiplying a matrix to the right and to the left by the permutation matrix $E=P_{\sigma}$ is to do the cycle $\sigma=(132)$ on the rows, and the inverse cycle $\sigma^{-1}=(123)$ on the columns. When doing this operation on a diagonal matrix, we obtain a non-diagonal matrix. In particular, the respective right hand sides of Eqs $(49)$ and $(50)$ cannot be the identity matrix. Moreover, since $W_i$ has order $2$, $\Delta$ has odd order and $W_i$ and $\Delta$ commute, their product cannot have order $2$. We conclude like announced. Then, we have
$$g(FUM)^3\in\lbrace (E^2FEF^{-1})^3,\,(FEF^{-1}E^2)^3,\,(E^2FE^2F^{-1}E^2)^3\rbrace$$
%Further, we claim that
%$$g(A)=F^2$$
Further, the $3$-Sylow $\mathcal{S}_3(F)$ is mapped to some $W_i$ conjugate of $\mathcal{S}_3(D)$. Now, an element of order $3$ of the center $Z(W_i\,\mathcal{S}_3(F)\,W_i^{-1})$ is a $W_i$ conjugate of an element of order $3$ belonging to $Z(\mathcal{S}_3(D))$, hence is $F^6$ or $F^{12}$ as shown earlier. Using some facts from before, it follows that
%Then, $g$ induces an isomorphism
%$$\bar{g}:3SylF/\mathcal{V}\;\longrightarrow\;3SylD/\mathcal{V}_D$$
%But $$3SylF/\mathcal{V}\;=\;\mathcal{S}_3(F)/\mathcal{N}$$
%and $$3SylD/\mathcal{V}_D\;=\;\mathcal{S}_3(D)/\mathcal{F}$$
%Then $g$ must map $\mathcal{N}$ onto $\mathcal{F}$ and $\mathcal{S}_3(F)$ onto $\mathcal{S}_3(D)$. In particular,
%$$\left|\begin{array}{l}g(H_3)\in\lbrace E,E^2\rbrace\\\\
%g(A)=F^2\end{array}\right.$$
%The fact that $g(A)=F^2$ is derived from before.
%Indeed, a study of the elements of order $3$ in the centers of the respective $3$-Sylows $\mathcal{S}_3(D)$ and $\mathcal{S}_3(F)$ had shown that
$$\left\lbrace\begin{array}{l}
g(A^3)=F^6\;\text{and}\;g(A^6)=F^{12}\\
\text{or}\\
g(A^3)=F^{12}\;\text{and}\;g(A^6)=F^6
\end{array}\right.$$
\newtheorem{Claim}{Claim}
\begin{Claim}
We have $F^6=\omega\,I_3$ and $F^{12}=\omega^2\,I_3$.
\end{Claim}
\noindent Hence it is not straightforward to conclude. We set for instance $g(A)=F^2$ like above.
%Since $\mathcal{N}$ is mapped onto $\mathcal{F}$ by $g$, the generator $A$ of $\mathcal{N}$ must be mapped to a matrix
%$$F^{2i}[F(F^{''})^{-1}]^{2j}$$
%for some integers $i$ and $j$. Then,
%$$g(A^3)=F^{6i}[F(F^{''})^{-1}]^{6j}\;,$$
%But $g(A^3)=F^{12}$ would imply that $g(A)=F^4$. Then, $g(A^2)=F^8\neq F^2$. So the second situation in the brace above is impossible and we must have $g(A)=F^2$.
%We are now in a position to show that $g$ maps $\mathcal{N}.\mathcal{V}$ onto the subgroup $N(18,1,1;2,1,1)$ of $D(18,1,1;2,1,1)$ generated by all the diagonal matrices. First,
Furthermore, we note that the imagae $g(H_3)$ cannot be in $\mathcal{F}$ since $\mathcal{F}$ commutes to $F^2$, but $H_3$ does not commute to $A$. Then it must be of the form
$$g(H_3)=W_i\Delta\,E^k\,W_i^{-1}$$
with $\Delta$ some diagonal matrix of $\mathcal{F}$, $k\in\lbrace 1,2\rbrace$, $i\in\lbrace 1,2,3,4\rbrace$ and $W_1=I_3$.

\noindent %Now the contradiction comes from the expression
We have the expression
$$g(C_6)=g(H_3)g(C_{18})g(H_3)^{-1}g(C_{18})^{8}$$

\noindent read out of the presentation given in Theorem $1$. Since $$C_{18}=A(FUM)^3,$$ the image $g(C_{18})$ is a diagonal matrix. The $E^k$ conjugate of a diagonal matrix is again a diagonal matrix. Therefore,
\begin{eqnarray}
g(H_3)g(C_{18})g(H_3)^{-1}&=&W_i\Delta\,E^k\,W_i^{-1}g(C_{18})W_iE^{-k}\Delta^{-1}W_i^{-1}\\
&=&E^k\,g(C_{18})E^{-k}
\end{eqnarray}
and so \begin{equation}
g(C_6)=E^k\,g(C_{18})E^{-k}g(C_{18})^{8}
\end{equation}
\\
%\noindent Moreover, since $C_6$ has order $6$, so must be the order of $g(C_6)$. Then,
%\begin{equation}
%g(C_6)=E^k\,g(C_{18})E^{-k}g(C_{18})^{-10}
%\end{equation}
By the presentation given in Theorem $1$, we also have
\begin{equation}
g(H_3)\,g(C_6)\,g(H_3)^{-1}=g(C_6)\,g(C_{18})^{-3}
\end{equation}
and so by similar arguments as before,
\begin{equation}
g(C_6)=E^k\,g(C_6)\,E^{-k}\,g(C_{18})^3
\end{equation}
By gathering Eqs $(53)$ and $(55)$, we now get
\begin{equation}
g(C_6)=E^{\bar{k}}\,g(C_{18})E^{-k}\,g(C_{18})^{8}E^{-k}g(C_{18})^3
\end{equation}
where $\bar{k}$ is the "conjugate" of $k$, namely $\bar{k}=1$ if $k=2$ and $\bar{k}=2$ if $k=1$. \\
With
$$\begin{array}{cc}\left|\begin{array}{l}
k=2\\
g((FUM)^3)=(E^2FEF^{-1})^3\\
g(H_1)=\overset{\sim}{B}\,E
\end{array}\right.&,\end{array}$$
all the relations of the presentation given in Theorem $1$ are verified on the images. Hence the Theorem.

\begin{Theorem}
The map
$$\begin{array}{ccc}
Fr(162\times 4)&\overset{g}{\longrightarrow}& D(18,1,1;2,1,1)\\
&&\\
H_3&\longmapsto& E^2\\
H_1&\longmapsto&\overset{\sim}{B}\,E\\
&&\\
C_{18}&\longmapsto&\begin{pmatrix}e^{-i7\pi/9}&&\\&e^{i2\pi/9}&\\&&e^{i5\pi/9}\end{pmatrix}\\
&&\\
C_6&\longmapsto&\begin{pmatrix}e^{i\pi/3}&&\\&-1&\\&&e^{i2\pi/3}\end{pmatrix}
\end{array}$$
where
$$\begin{array}{ccc}
g(C_{18})&=&F^2\,(E^2FEF^{-1})^3\\
g(C_6)&=&E\,g(C_{18})E^{-2}\,g(C_{18})^{8}E^{-2}\,g(C_{18})^3\\
&&
\end{array}$$
defines an isomorphism of groups between $Fr(162\times 4)$ and $D(18,1,1;2,1,1)$.
\end{Theorem}
\section{Epilogue}
We have just shown that the groups $D(18,1,1;2,1,1)$ and $Fr(162\times 4)$ are isomorphic groups and both groups are isomorphic to a semi-direct product $(\mathbb{Z}_6\times\mathbb{Z}_{18})\rtimes S_3$. In this part, we show that the only $SU(3)$-subgroups from the extended version of the $1916$ classification the Freedman group may be isomorphic to are of $D$-type and we provide some instances of $D$-groups the Freedman group is isomorphic to.
At first sight, we can already rule out the series $\Delta(3n^2)$ and $\Delta(6n^2)$ because
$$648=6\times 108=3\times 216$$
but neither of $108$ or $216$ is a square. Further, for orders considerations, the only candidate in the series $T_n$ would be $T_{216}$ but the series exists only for special values of the integer $n$ (see \cite{FF1}) and $216$ being even is not a product of primes of the form $6k+1$ with $k\geq 1$ an integer. The only candidate among the exceptional groups or direct products of $\mathbb{Z}_3$ by an exceptional group would in turn be $\Sigma(216\times 3)$, still for orders purposes. It is a result from P.O. Ludl thesis that this group is generated by only two matrices which we recall below.
$$\Sigma(216\times 3)=\bigg\langle D=\begin{pmatrix} \varepsilon&&\\&\varepsilon&\\&&\varepsilon\omega\end{pmatrix}\,,\,V=\frac{1}{i\sqrt{3}}\begin{pmatrix}
1&1&1\\1&\omega&\omega^2\\1&\omega^2&\omega\end{pmatrix}\bigg\rangle$$
where $$\varepsilon=e^{\frac{4i\pi}{9}}\;\text{and}\;\omega=e^{\frac{2i\pi}{3}}$$
However, we have the following result.
\begin{Theorem}
The Freedman group $Fr(162\times 4)$ is not isomorphic to the exceptional group $\Sigma(216\times 3)$.
\end{Theorem}
\textsc{Proof.} We have the following principal series, where the numbers at the top denote the degrees of the group extensions.
$$\negthickspace\negthickspace\negthickspace\negthickspace\negthickspace\negthickspace\negthickspace\negthickspace
\negthickspace\negthickspace\negthickspace\negthickspace\negthickspace\negthickspace\negthickspace\negthickspace
\negthickspace\negthickspace\negthickspace\negthickspace\negthickspace\negthickspace\negthickspace\negthickspace
\begin{array}{ccccccccccccccc}
&3&&3&&3&&2&&2&&2&&3&\\
\lbrace e\rbrace&\lhd&\mathbb{Z}_3&\lhd&\mathbb{Z}_3\times\mathbb{Z}_3&\lhd&\Delta(27)&\lhd&\Delta(54)&\lhd&\Sigma(36\times 3)&\lhd&\Sigma(72\times 3)&\lhd&\Sigma(216\times 3)\\\\
&&&&&&&&&&&&&&\\
&3&&3&&3&&2&&2&&3&&2&\\
\lbrace e\rbrace&\lhd& <A^3>&\lhd&<A>&\lhd&\mathcal{N}&\lhd&\mathcal{N}.<(FUM)^3>&\lhd&\mathcal{N}.\mathcal{V}&\lhd& 3SylF&\lhd& Fr(162\times 4)
\end{array}$$
%$$\begin{array}{l}
%&3&&9&&2&&4&&3&\\
%\lbrace e\rbrace\;\lhd\;\mathbb{Z}_3\;\lhd\;\mathbb{Z}_3\times\mathbb{Z}_3\;\lhd\;\Delta(27)\;\lhd\;\Delta(54)\;\lhd\;\Sigma(36\times 3)\;\lhd\;\Sigma(72\times 3)\;\lhd\;\Sigma(216\times 3)\\\\
%&&&&&&&&&&&&\\
%&9&&3&&4&&3&&2&\\
%\lbrace e\rbrace\;\lhd\; <A^3>\;\lhd\;<A>\;\lhd\;\;\mathcal{N}\;\lhd\;\mathcal{N}.<(FUM)^3>\;\lhd\;\mathcal{N}.\mathcal{V}\;\lhd\; 3SylF\;\lhd\; Fr(162\times 4)
%\end{array}$$

\noindent The first one is a result of \cite{GL} where the principal series of the exceptional $SU(3)$ finite subgroups are determined.
The second one is derived from the present paper. It is also shown in \cite{GL} that $\Sigma(72\times 3)/\Delta(27)$ is isomorphic to the quaternion group $Q_8$.
\begin{Lemma}
Suppose $Fr(162\times 4)$ contains a proper subgroup $H$ such that $H$ has a quotient $H/K$ with $H/K\simeq Q_8$. Then, all the $2$-Sylow subgroups of $Fr(162\times 4)$ are isomorphic to $Q_8$.
\end{Lemma}
\textsc{Proof.} If such subgroups $H$ and $K$ exist, then $2^3$ divides $|H|$, and so
$$|H|=2^3.\,3^k$$
for some integer $k\in\lbrace 0,1,2,3\rbrace$ and $|K|=3^k$. First, if $k=0$, the result clearly holds. Assume now $k\in\lbrace 1,2,3\rbrace$. Since $K$ is normal in $H$, then $K$ must be the unique $3$-Sylow subgroup $S_3$ of $H$. Let $S_2$ be a $2$-Sylow subgroup of $H$. Since $S_2\cap S_3=\lbrace e\rbrace$, we then have
$$S_3.S_2=H\simeq S_3\rtimes S_2$$
It follows that $$S_2\simeq H/S_3\simeq Q_8$$
Thus, $S_2$ is isomorphic to $Q_8$. Further, $S_2$ is also a $2$-Sylow subgroup of $Fr(162\times 4)$ and all the $2$-Sylows of $Fr(162\times 4)$ are isomorphic since they are all conjugate. We deduce that any $2$-Sylow subgroup of $Fr(162\times 4)$ must then be isomorphic to the quaternion group $Q_8$.\hfill$\square$
\begin{Lemma}
$\mathcal{V}\rtimes <H_1>$ is a $2$-Sylow subgroup of $Fr(162\times 4)$ which is not isomorphic to $Q_8$.
\end{Lemma}
\textsc{Proof.} Recall $\mathcal{V}$ is a normal subgroup of $Fr(162\times 4)$. The semi-direct product from the statement of the Lemma has the right cardinality, hence is a $2$-Sylow subgroup of $Fr(162\times 4)$. One of the characteristics of $Q_8$ is that all its subgroups are normal. However,
$$V_2.\,H_1\,.V_2^{-1}\not\in\,<H_1>,$$
so that $<H_1>$ is a subgroup of $\mathcal{V}\rtimes <H_1>$ which is not normal. This shows the Lemma. \\
Now, Lemma $9$ and Lemma $10$ imply Theorem $9$. \hfill $\square$\\\\
%We see that the number of terms in the respective two series differs by one, thus by Jordan-H\"older Theorem the two groups $Fr(162\times 4)$ and $\Sigma(216\times 3)$ cannot be isomorphic.\hfill $\square$
%We use further Jordan-H\"older Theorem to
 We now study the $SU(3)$-subgroups from the series $(D)$. It is known only since recently $(2011)$ in the work of \cite{PO2} that an $SU(3)$ finite subgroup from the series $(D)$ has the general structure
$$\mathcal{A}\rtimes S_3\simeq \mathbb{Z}_m\times\mathbb{Z}_p\rtimes S_3$$
with $\mathcal{A}$ the normal subgroup of all the diagonal matrices of the given $D$-group.
Moreover, $m$ divides $p$ by Theorem
$2.1$ of \cite{PO2} which states that any finite abelian subgroup of $SU(3)$ is of the form $\mathbb{Z}_r\times\mathbb{Z}_s$ with $r$ divides $s$. Furthermore, we prove the following Lemma.
\begin{Lemma}
The Freedman group $Fr(162\times 4)$ does not contain any element of order $24$, $27$, $54$, $72$ or $108$.
\end{Lemma}
\textsc{Proof.} Suppose it does contain an element of order $k$ with $$k\in\lbrace 24,27,54,72,108\rbrace$$
First, there does not exist any element of such order inside the direct product $\mathbb{Z}_3\times\mathbb{Z}_9\times V_4$. Indeed, if $d_1|3$, $d_2|9$ and $d_3|2$, then none of
$$24,27,54,72,108$$
does divide
$$lcm(d_1,d_2,d_3)$$
%$$lcm(d_1,d_2,d_3)\not\in\lbrace 24,27,54,72,108\rbrace$$
So an element of order $k$ must be a product
$$N\,H_i$$
with $N$ in the direct product and $H_i$ in the symmetric group.
Moreover, by the way the semi-direct product is defined, we must have
$$H_i^k=I_3$$
This implies $i\in\lbrace 0,1,2\rbrace$ if $k$ is even and $i\in\lbrace 3,4\rbrace$ if $k$ is odd.\\
%Further, still by the way the semi-direct product is defined, we must have
%$$H_i=NH_iN^{k-1}$$
A generic form for $N$ is
$$\begin{array}{l}N=A^sB^tV_r\\\end{array}$$
Let us first deal with the case $k=27$. We found out using Mathematica that all the products $NH_3$ or $NH_4$ have order $3$. So there is no element of order $27$ in the Freedman group. \\
%We have $$N^{27}=V_r^{27}=V_r$$
%It follows that
%$$N=H_iV_r^{-1}N\,H_i^{-1}\;,$$
%and so,
%$N$ and $V_rN$ being conjugate have the same order. In other words, $N$ and $A^sB^t$ have the same order. This implies that $N$ has an odd order. Consequently, $N=A^sB^t$. Further, by inspection, the only elements $A^sB^t$ which commute to $H_i$ are $A^3$ and $A^6$. Moreover, all the products $A^3H_i$ and $A^6H_i$ have order $3$ instead of $27$. We conclude that there is no element of order $27$ in $Fr(162\times 4)$.
%None of the products $A^sB^tH_i$ with $i\in\lbrace 3,4\rbrace$ and $r\in\lbrace 2,3,4\rbrace$ has order $27$. These products all have order $3$. \\\\
Since $54=27\times 2$ and $108=54\times 2$, there are no elements of order $54$ or $108$ either since their respective squares would then have order $27$ and $54$. %\\If the Freedman group contained an element of order $108$, its square would have order $54$, which, like just seen, is impossible.
\\Finally, we see with Mathematica that the possible orders for the products $NH_i$ with $N\in\mathcal{N}\mathcal{V}$ and $i\in\lbrace 0,1,2\rbrace$ are not among $24$ or $72$. \hfill $\square$\\\\
Since \begin{eqnarray*} 108&=&2\times 54\\
&=&3\times 36\\
&=&4\times 27\\
&=&6\times 18\\
&=&9\times 12\end{eqnarray*}
we deduce that the only candidates from the series $(D)$ of order $648$ are of the form
$$\begin{array}{l}\mathbb{Z}_3\times\mathbb{Z}_9\times\mathbb{Z}_2\times\mathbb{Z}_2\rtimes S_3\\
\mathbb{Z}_3\times\mathbb{Z}_9\times\mathbb{Z}_4\rtimes S_3
\end{array}$$

We now prove the following Theorem.
\begin{Theorem}\hfill\\
$i)$ Assume without loss of generality that $gcd(n,a,b)=1$. %If the Freedman group is isomorphic to a $D$-group $D(n,a,b;d,r,s)$ with
%$$\left|\begin{array}{l}gcd(n,a,b)=1\\ r\neq 0\;\;\text{ when $d$ is odd}\\ r\neq \frac{d}{2}\;\;\text{ when $n$ is even,}\end{array}\right.$$
A necessary condition for the Freedman group to be isomorphic to the $D$-group $D(n,a,b;d,r,s)$ is that
$$n\in\lbrace 2,3,6,9,18\rbrace$$
$ii)$ The Freedman group is isomorphic to many $SU(3)$ subgroups from the series $(D)$. For instance, it is isomorphic to all the identical groups
$$\begin{array}{l}D(9,1,1;2,0,1),\;D(9,1,1;2,0,0),\\D(18,1,1;2,1,1),\;D(18,1,1;2,0,0),\;D(18,1,1;2,0,1)
\end{array}$$
\end{Theorem}
\textsc{Proof.}
%D(n,a,b;d,r,s)\;\text{contains a subgroup which is isomorphic to}\;D(9,1,1;2,1,1)\\
%D(n,a,b;d,r,s)\;\text{is not isomorphic to}\;D(18,1,1;2,1,1)
%\end{array}$$
Notice that
$$C(n,a,b)\,\subseteq\,D(n,a,b;d,r,s)$$
As shown in \cite{GL}, the structure of an $SU(3)$ subgroup from the series $(C)$ is
$$\mathcal{A}\rtimes \mathbb{Z}_3\simeq(\mathbb{Z}_m\times\mathbb{Z}_p)\rtimes\mathbb{Z}_3$$
with $\mathcal{A}$ the normal subgroup of all the diagonal matrices of $C(n,a,b)$, the integer $p$ the maximal order of the diagonal matrices and $m$ the divisor of $p$ provided in \cite{PO2}.\\
In particular, in a $C$-group, all the elements of order $2$ commute. Thus, the inclusion above is strict.
If $|C(n,a,b)|$ divides strictly $3^4.2^3$, this implies that
$$|C(n,a,b)|\leq 324$$
%In \cite{PO3}, P.O. Ludl provides a table which contains a list of all the finite $(C)$-subgroups of $SU(3)$ of order less than $512$. A $(C)$-group from this table whose order is less than or equal to $324$ and divides $648$ can only be one of the following, up to isomorphism.
%$$\begin{array}{l}C(l,0,0)\\l>0\end{array},\;\Delta(12), \;\Delta(27),\;\Delta(108),\;C(9,1,1),\;C(18,1,1)$$
%\newtheorem{Proposition}{Proposition}
%\begin{Proposition}\textbf{The following statement is not up to isomorphism}.\\
%Up to multiplication of the three defining integers by a common non-zero integer, these are the only $C$-groups which may occur inside a $D$-group isomorphic to the Freedman group.
%\end{Proposition}
%\textsc{Proof.}
Further, let $g=gcd(n,a,b)$. Notice that
$$C\bigg(\frac{n}{g},\frac{a}{g},\frac{b}{g}\bigg)=C(n,a,b)$$
Thus, without loss of generality, we may assume that $gcd(n,a,b)=1$ and so the order of $F(n,a,b)$ is $n$. From before, a diagonal matrix of a $D$-group isomorphic to the Freedman group has order among
$$2,3,4,6,9,12,18,36$$
And so,
$$\begin{array}{cc}n\in\lbrace 2,3,6,9,18\rbrace&(\star)\\\qquad\text{or}\qquad&\\$$
n\in\lbrace 4,12,36\rbrace& (\star\star)\end{array}$$

\noindent We claim that the second row of values are to exclude. Our proof is based on the following set of propositions.
\newtheorem{Proposition}{Proposition}
\begin{Proposition}
We have
$$C(k,\bar{a},\bar{b})\subseteq C(2k,a,b)$$
with $$\bar{a}=\begin{cases} a&\text{if $0\leq a\leq k-1$}\\a-k&\text{if $k\leq a\leq 2k-1$}\end{cases}$$
We also have
$$C(k,\hat{a},\hat{b})\subseteq C(3k,a,b)$$
with $$\hat{a}=\begin{cases} a&\text{if $0\leq a\leq k-1$}\\a-k&\text{if $k\leq a\leq 2k-1$}\\
a-2k&\text{if $2k\leq a\leq 3k-1$}\end{cases}$$
\end{Proposition}
\noindent\textsc{Proof.} Straightforward.
\begin{Proposition}
$$C(2,0,1)=C(2,1,0)=C(2,1,1)$$
Moreover, these groups contain a subgroup which is isomorphic to the Klein group $V_4$.
\end{Proposition}
\textsc{Proof.} It suffices to notice that the matrices $F(2,0,1)$, $F(2,1,0)$ and $F(2,1,1)$ are all conjugate under $E$ or $E^2$.
\begin{Proposition} Suppose $gcd(4,a,b)=1$.
A $D$-group isomorphic to the Freedman group cannot contain a $C$-group of type $C(4,a,b)$.
\end{Proposition}
\textsc{Proof.} Suppose it does. Then it also contains a $C$-group of type $C(2,\bar{a},\bar{b})$ by Proposition $1$. Then, by Proposition $2$, the group $C(4,a,b)$ contains a Klein group formed of diagonal matrices. Since $gcd(4,a,b)=1$, the group $C(4,a,b)$ must also contain a diagonal matrix of order $4$. Then, the $D$-group containing $C(4,a,b)$ would contain a direct product $\mathbb{Z}_4\times\mathbb{Z}_2\times\mathbb{Z}_2$ of diagonal matrices. However, we know from before that the maximal power of two dividing the order of the subgroup of diagonal matrices is $4$, thus a contradiction.
\begin{Proposition} Suppose $gcd(12,a,b)=1$.
A $D$-group isomorphic to the Freedman group cannot contain a $C$-group of type $C(12,a,b)$.
\end{Proposition}
\textsc{Proof.} Follows from Propositions $1$ and $3$.
\begin{Proposition} Suppose $gcd(36,a,b)=1$.
A $D$-group isomorphic to the Freedman group cannot contain a $C$-group of type $C(36,a,b)$.
\end{Proposition}
\textsc{Proof.} Follows from Propositions $1$ and $4$.\\\\
This finishes the proof of point $(i)$ in the Theorem. Our study of which $D$-groups arise as the Freedman group is far from being complete at this point. However, point $(i)$ and its proof provide in fact many more informations than those already disclosed about the integers $n$, $a$ and $b$. Let us give an example.
\newtheorem{Fact}{Fact}
\begin{Fact}Suppose $s$ is an integer such that $(3,s)=1$. \\
The group $C(9,3,s)$ never arises as a subgroup of the Freedman group. \\
Consequently also, nor does the group $C(18,3,s)$ arise as a subgroup of the Freedman group.
\end{Fact}
\noindent\textsc{Proof.} The second point follows from the first one by applying Proposition $1$.
%We will show that when $n=18$ we could for instance never have $a=7$ and $b=5$.
We will show that when $(3,s)=1$, the structure of $C(9,3,s)$ is
$$\mathbb{Z}_{9}\times\mathbb{Z}_9\rtimes \mathbb{Z}_3$$
Then we have $3^5$ divides the order of $C(9,3,s)$, which prevents this group from being isomorphic to a subgroup of the Freedman group. \\Since $gcd(9,3,s)=1$, the maximal order of a diagonal matrix of $C(9,3,s)$ is $9$. Hence
$$C(9,3,s)\simeq \mathbb{Z}_{9}\times\mathbb{Z}_p\rtimes \mathbb{Z}_3,$$
with $p\in\lbrace 3,9\rbrace$.
Following \cite{PO2}, the integer $p$ must satisfy to
$$\left\lbrace\begin{array}{cccc}
9&|&p(s-3t)&1\leq t\leq \frac{9}{p}-1,\;\;\text{and $p$ the smallest}\\
9&|&p(3+s(t+1))&
\end{array}\right.$$
If $p=3$, then we must have $3|s-3t$ which implies $3|s$, impossible. Thus, we rather have $p=9$. \hfill $\square$\\

\noindent We provide below more inclusions of $C$-groups inside $D$-groups which impose more restrictions on the choice for the integers $d$ and $r$. First, we deal with the case when $d$ is even.

\begin{Lemma} Assume $d$ is even. Then, the following inclusions hold.
$$\left\lbrace\begin{array}{ccccc}
C(d,2r,\frac{d-2r}{2})&\subset &D(n,a,b;d,r,s)&\text{if}& r\leq\frac{d}{2}-1\\
&&&&\\
C(d,2r-d,\frac{3d-2r}{2})&\subset& D(n,a,b;d,r,s)&\text{if}& r\geq \frac{d}{2}+1\\
&&&&\\
C(d,0,0)&\subset&D(n,a,b;d,\frac{d}{2},s)&\text{if}&r=\frac{d}{2}
\end{array}\right.$$

\end{Lemma}
\textsc{Proof.} Read from
$$\overset{\sim}{G}(d,r,s)^2=\begin{pmatrix}
e^{2i2r\pi/d}&0&0\\0&-e^{-2i\pi r/d}&0\\0&0&-e^{-2i\pi r/d}
\end{pmatrix}\\$$
\hfill $\square$\\
\noindent The case $d$ is odd requires a bit more effort. Instead, we must consider the fourth power of the generator $\overset{\sim}{G}(d,r,s)$. We have
$$\overset{\sim}{G}(d,r,s)^4=\begin{pmatrix} e^{\frac{2i\pi}{d}4r}&&\\&e^{\frac{2i\pi}{d}(-2r)}&\\&&e^{\frac{2i\pi}{d}(-2r)}\end{pmatrix}$$
From there, we deduce the following lemma.
\begin{Lemma} Suppose $d$ is odd. The following inclusions hold.
$$\left\lbrace\begin{array}{ccccc}
C(d,4r-3d,2d-2r)&\subset&D(n,a,b;d,r,s)&\text{if}&d\geq 5\;\;\&\;\;r\geq\lfloor\frac{3d}{4}\rfloor +1 \\
&&&&\\
C(d,d-3,\frac{d+3}{2})&\subset& D(n,a,b;d,\lfloor\frac{3d}{4}\rfloor,s)&\text{if}&d\geq 5\;\;\&\;\;r=\lfloor\frac{3d}{4}\rfloor\;\;\&\;\;d\equiv\;1\;mod\;4\\
&&&&\\
C(d,d-1,\frac{d+1}{2})&\subset& D(n,a,b;d,\lfloor\frac{3d}{4}\rfloor,s)&\text{if}&d\geq 3\;\;\&\;\;r=\lfloor\frac{3d}{4}\rfloor\;\;\&\;\;d\equiv\;3\;mod\;4\\
&&&&\\
C(d,4r-2d,2d-2r)&\subset& D(n,a,b;d,r,s)&\text{if}&d\geq 7\;\;\&\;\;\lfloor\frac{d}{2}\rfloor +1\leq r\leq\lfloor\frac{3d}{4}\rfloor-1\\
&&&&\\
C(d,d-2,1)&\subset&D(n,a,b;d,\lfloor\frac{d}{2}\rfloor,s)&\text{if}&d\geq 3\;\;\&\;\;r=\lfloor\frac{d}{2}\rfloor\\
&&&&\\
C(d,4r-d,d-2r)&\subset&D(n,a,b;d,r,s)&\text{if}&d\geq 7\;\;\&\;\;\lfloor\frac{d}{4}\rfloor+1\leq r\leq\lfloor\frac{d}{2}\rfloor-1\\
&&&&\\
C(d,d-1,\frac{d+1}{2})&\subset&D(n,a,b;d,\lfloor\frac{d}{4}\rfloor,s)&\text{if}&d\geq 3\;\;\&\;\;r=\lfloor\frac{d}{4}\rfloor\;\;\&\;\;d\equiv 1\; mod\; 4\\
&&&&\\
C(d,d-3,\frac{d+3}{2})&\subset& D(n,a,b;d,\lfloor\frac{d}{4}\rfloor,s)&\text{if}&d\geq 5\;\;\&\;\;r=\lfloor\frac{d}{4}\rfloor\;\;\&\;\;d\equiv 3\;mod\;4\\
&&&&\\
C(d,4r,d-2r)&\subset&D(n,a,b;d,r,s)&\text{if}&d\geq 13\;\;\&\;\;1< r\leq \lfloor\frac{d}{4}\rfloor -1\\
&&&&\\
C(d,4,d-2)&\subset&D(n,a,b;d,1,s)&\text{if}&d\geq 5\;\;\&\;\;r=1\\
&&&&\\
C(3,1,1)&\subset&D(n,a,b;3,1,s)&\text{if}&d=3\;\;\&\;\;r=1\\
&&&&\\
C(d,0,0)&\subset&D(n,a,b;d,0,s)&\text{if}&r=0
\end{array}\right.$$
where $\lfloor x\rfloor$ denotes the integer part of $x$.\\
\end{Lemma}
%This ends the proof of Theorem $10$. Our proof also implies the following result.

\noindent Knowing from Ludl's work in \cite{PO2} the structure of the $C$-groups, we obtain this way many informations on the structure of the $D$-groups.\\
We next focus our attention on some $D$-groups $D(n,a,b;d,r,s)$ with $n\in\lbrace 9,18\rbrace$, $(a,b)=(1,1)$ and $d=2$. We have
$$\overset{\sim}{G}(2,0,s)=\begin{pmatrix}1&0&0\\0&0&(-1)^s\\0&(-1)^{s+1}&0\end{pmatrix}$$
We thus see that the two groups $D(18,1,1;2,0,0)$ and $D(18,1,1;2,0,1)$ are identical since
$$\overset{\sim}{G}(2,0,0)=\overset{\sim}{G}(2,0,1)^{-1}$$
Moreover, we notice that
$$\overset{\sim}{B}\,W_4=\overset{\sim}{G}(2,0,1)$$
Thus we have
$$D(18,1,1;2,0,1)\subseteq D(18,1,1;2,1,1)$$
The latter inclusion is in fact an equality since we have
$$F(18,1,1)^9=W_4$$
So the three groups $D(18,1,1;2,1,1)$, $D(18,1,1;2,0,1)$ and $D(18,1,1;2,0,0)$ are identical. %In order to conclude this case, it suffices to prove the following fact.
\begin{Fact}
For the smaller groups, we have
$$D(9,1,1;2,1,1)\subset D(9,1,1;2,0,1)=D(9,1,1;2,0,0)=D(18,1,1;2,1,1)$$
\end{Fact}
\noindent \textsc{Proof.} We already know that the first equality holds. Next, notice that
$$(\overset{\sim}{G}(2,0,1)E)^2=W_4$$
Hence $$W_4\in D(9,1,1;,2,0,1)$$
Then $\overset{\sim}{B}\in D(9,1,1;2,0,1)$ and so
$$D(9,1,1;2,1,1)\subseteq D(9,1,1;2,0,1)\subseteq D(18,1,1;2,1,1)$$
Because in $D(9,1,1;2,1,1)$, there is no diagonal matrix of even order, we know that $W_4$ cannot belong to $D(9,1,1;2,1,1)$. Hence the first inclusion above is strict. \\Then, either $|D(9,1,1;2,0,1)|=324$ or $D(9,1,1;2,0,1)=D(18,1,1;2,1,1)$. To decide, it suffices to notice that $W_4$ belongs to $D(9,1,1;2,0,1)$ implies that the whole Klein group $\mathcal{V}_D$ is contained in $D(9,1,1;2,0,1)$. Indeed, the matrices $W_2$ and $W_3$ are the respective $E$ and $E^2$ conjugates of the matrix $W_4$. \\By \cite{PO2}, the D-group $D(9,1,1;2,0,1)$ has the structure
$$\mathcal{A}\rtimes S_3$$
with $\mathcal{A}$ the normal subgroup of all the diagonal matrices of $D(9,1,1;2,0,1)$. \\So,
$$4\,\Big{|}\;|\mathcal{A}|$$
Then,
$$24\,\Big{|}\;|D(9,1,1;2,0,1)|$$
This prevents to have
$$|D(9,1,1;2,0,1)|=324$$
And so,
$$D(9,1,1;2,0,1)=D(18,1,1;2,1,1)$$
This ends the proof of Theorem $10$. \\

 %Also, we finish this part on $D$-groups by giving a conjecture.
%\newtheorem{Conjecture}{Conjecture}
%\begin{Conjecture}
%The semi-direct product $\mathbb{Z}_6\times\mathbb{Z}_{18}\rtimes S_3$ arises only as five identical groups in the series $(D)$, namely as
%$$D(9,1,1;2,0,1),\;D(9,1,1;2,0,0),\;D(18,1,1;2,1,1),\; D(18,1,1;2,0,0),\;D(18,1,1;2,0,1)$$
%Thus,
%Out of the $36$ non-isomorphic semi-direct products $\mathbb{Z}_6\times\mathbb{Z}_{18}\rtimes S_3$, a single one of them is a $(D)$-group and is the Freedman group.
%\end{Conjecture}
%\textsc{Proof.}\\\\ By earlier considerations, a $(D)$-group of that order must contain a diagonal matrix of order $9$ and must be so that the $(C)$-groups which it contains are among $C(2,0,1)$, $C(3,0,1)$, $C(6,0,1)$, $C(9,1,1)$, $C(18,1,1)$. Only the five identical $(D)$-groups from the statement of the Theorem do satisfy to these criteria. \\
%As for the last statement in the Theorem,
%We wrote a program in GAP which uses the small groups library \cite{GAP} and shows that there are exactly $36$ non-isomorphic semi-direct products
%$$\mathbb{Z}_6\times\mathbb{Z}_{18}\rtimes S_3$$
\indent \textsc{Remarks}.\\\\
i) The order of the direct product $\mathbb{Z}_6\times\mathbb{Z}_{18}$ is the number of distinct semi-direct products $\mathbb{Z}_6\times\mathbb{Z}_{18}\rtimes S_3$ that are all isomorphic to the Freedman group.\\
ii) The order of the automorphism group $Aut(\mathbb{Z}_6\times\mathbb{Z}_{18})$ is the order of the Freedman group. \\\\
Point i) was obtained by writing a program in GAP which lists all the semi-direct products $\mathbb{Z}_6\times\mathbb{Z}_{18}\rtimes S_3$ from the list of all the homomorphisms from $S_3$ to $Aut(\mathbb{Z}_6\times\mathbb{Z}_{18})$ and then returns the list of their GAP ID. There are $36$ distinct GAP ID's in this list, that is there are $36$ non-isomorphic semi-direct products $\mathbb{Z}_6\times\mathbb{Z}_{18}\rtimes S_3$. We counted the number of occurrences in this list of the GAP ID
$[468,259]$ (see forthcoming (63)) and found the number $108$. \\
%\begin{Lemma} Suppose $d$ is even. Then we have,
%$$\left\lbrace\begin{array}{ccc}
%C(d,2r,\frac{d-2r}{2},\frac{d-2r}{2}\subset D(n,a,b

%\end{Lemma}
%We now investigate the series $(C)$.
%Though we did not need it earlier and won't need this result in the forthcoming discussions, it is worth mentioning the following result.

%We then deduce the following statement.
We now investigate the series $(C)$. We have the following statement whose proof is straightforward.
\begin{Theorem}
The Freedman group is not isomorphic to any $SU(3)$ finite subgroup from the series $(C)$.
\end{Theorem}
\textsc{Proof.} %First we recall Theorem $2.1$ from \cite{PO2} which states that any finite abelian subgroup of $SU(3)$ is of the form $\mathbb{Z}_n\times\mathbb{Z}_m$ with $n$ divides $m$. Another
As already mentioned before, a result of \cite{PO2} states that a $(C)$ group is of the form $\mathbb{Z}_n\times\mathbb{Z}_m\rtimes\mathbb{Z}_3$. Hence, in such a group all the elements of even order commute. In a $(D)$-group however, the elements of order two of the symmetric group do not commute.\\\\  %Gathering both results, we must have $n$ divides $m$ in the expression for the latter group. By using our Lemma $14$, it appears that the only candidate would then be $\mathbb{Z}_6\times\mathbb{Z}_{36}\rtimes\mathbb{Z}_3$. Now simply notice that any element of order $4$ in this group commutes to any element of order $2$ in the group. But in the Freedman group $V_1H_0$ is an element of order $4$ which does not commute to  \\\\
Finally, the Freedman group is not isomorphic to a group of matrices of the shape
$$\begin{array}{l}\begin{array}{cc}\begin{pmatrix}
det(A^{\dagger})& 0_{1,2}\\
0_{2,1}&A
\end{pmatrix}&,\end{array}\\\\\end{array}$$
with $A$ a unitary matrix of size $2$. So it does not belong to the "B type" as it is referred to in the classification of \cite{PO2}. Indeed, we have the following statement.
\begin{Theorem}
The Freedman group is not isomorphic to any $SU(3)$ finite subgroup of the type $B$. In other words, the Freedman group is not isomorphic to a finite subgroup of $U(2)$.
\end{Theorem}
However, the Freedman group contains subgroups of type $(B)$ of various orders. It is worth mentioning one which has order $54$.
\begin{Theorem}
The Freedman group contains a $(B)$-subgroup which is isomorphic to $\mathbb{Z}_9\times S_3$.
\end{Theorem}
%Along the way, we also show the following statement.
%\begin{Theorem}
%The Freedman group is not an extension of $\Delta(54)$.
%\end{Theorem}
%Thus recovering the fact that $Fr(162\times 4)$ and $\Sigma(216\times 3)$ are not isomorphic groups. \\
%It is also worth mentioning at this point that $Fr(162\times 4)$ is contained in the $\Delta$ subseries of the series  %$(D)$ by the main Theorem $7$ and the remark which follows.
%\newtheorem{Remark}{Remark}
%\begin{Remark}
%In $2009$, R. Zwicky and T. Fischbacher showed in their paper \cite{FZ} that every $(D)$-group is a subgroup of $\Delta(6\,n^2)$ for a suitable $n$. Then the Freedman group also arises as a subgroup of some $\Delta(6n^2)$.
%\end{Remark}

\textsc{Proof.} We first begin with the latter Theorem. Our proof here relies on a table which can be found by clicking on a link created by the same authors in the $2011$ paper \cite{PW} by Parattu and Wingerter. This table shows in particular whether groups of GAP ID $[54,n]$ are subgroups of $SU(3)$ or not and whether they are subgroups of $U(2)$ or not. We summarize below the status with respect to $U(2)$ for those groups belonging to $SU(3)$. The informations below have been copied out of Parattu and Wingerter's table. We added an extra row and extra column to the table for the sake of $(\mathbb{Z}_9\times\mathbb{Z}_3)\rtimes\mathbb{Z}_2$, which drew our attention while writing this proof. \\
\begin{center}
\begin{tabular}{|c|c|c|c|}
\hline
\text{GAP ID}&\text{Group}&SU(3)&U(2)\\
\hline
[54,1]&$D_{27}$&YES&YES\\
\hline
[54,3]&$\mathbb{Z}_3\times D_9$&YES&YES\\
\hline
[54,4]&$\mathbb{Z}_9\times S_3$&YES&YES\\
\hline
[54,7]&$(\mathbb{Z}_9\times\mathbb{Z}_3)\rtimes\mathbb{Z}_2$&NO&NO\\
\hline
[54,8]&$((\mathbb{Z}_3\times\mathbb{Z}_3)\rtimes\mathbb{Z}_3)\rtimes\mathbb{Z}_2$&YES&NO\\
\hline
\end{tabular}
\end{center}
Our study of the Freedman group shows that it contains a subgroup $$S=\mathcal{N}\rtimes <H_i>$$ with $i\in\lbrace 0,1,2\rbrace$, which is isomorphic to a semi-direct product
$$(\mathbb{Z}_9\times\mathbb{Z}_3)\rtimes\mathbb{Z}_2$$
We notice the semi-direct product $(\mathbb{Z}_9\times\mathbb{Z}_3)\rtimes\mathbb{Z}_2$ which is recorded in Parattu and Wingerter's table does not arise as a subgroup of $SU(3)$.
This means the two semi-direct products are not isomorphic. By Lemma $11$, the Freedman group does not contain any element of order $27$. Thus, $S$ cannot be isomorphic to $D_{27}$ either. By \cite{PO2}, the group of GAP ID $[54,8]$ is $\Delta(54)$. In order to rule it out, it suffices to notice that a group of order $54$ has a unique $3$-Sylow subgroup. In the case of $S$, this unique $3$-Sylow subgroup is the group $\mathcal{N}$ of order $3^3$. Recall that the groups of order $3^3$ are up to isomorphism

$$\begin{array}{cc}\begin{array}{l}\text{The two non-abelian groups}\\
\begin{array}{cc}
(\mathbb{Z}_3\times\mathbb{Z}_3)\rtimes\mathbb{Z}_3&
\mathbb{Z}_9\rtimes\mathbb{Z}_3\end{array}\end{array}&
\begin{array}{l}\text{The three abelian groups}\\
\begin{array}{ccc}
\mathbb{Z}_3\times\mathbb{Z}_3\times\mathbb{Z}_3&
\mathbb{Z}_9\times\mathbb{Z}_3&
\mathbb{Z}_{27}
\end{array}\end{array}\end{array}$$
Now, the group $S$ cannot be isomorphic to $((\mathbb{Z}_3\times\mathbb{Z}_3)\rtimes\mathbb{Z}_3)\rtimes\mathbb{Z}_2$, else it would contain a non-abelian $3$-Sylow subgroup. %or one which is isomorphic to $\mathbb{Z}_3\times\mathbb{Z}_3\times\mathbb{Z}_3$.
However, its unique $3$-Sylow subgroup is $\mathcal{N}$ which is abelian. %and
%$$\mathcal{N}\not\simeq\mathbb{Z}_3\times\mathbb{Z}_3\times\mathbb{Z}_3$$
Note that we also recover the fact that $S$ is not isomorphic to
$D_{27}=\mathbb{Z}_{27}\rtimes\mathbb{Z}_2$, by a similar argument. %And of course, $S$ is not isomorphic to $\mathbb{Z}_3\times\mathbb{Z}_3\times\mathbb{Z}_3\rtimes\mathbb{Z}_2$ since
%$$\mathcal{N}\not\simeq\mathbb{Z}_3\times\mathbb{Z}_3\times\mathbb{Z}_3$$
By combining all the preceding results and the content of the table above, the group $S$ must then be isomorphic to one of the two groups
$$\begin{array}{l}
\mathbb{Z}_3\times D_9\\\\
\mathbb{Z}_9\times S_3
\end{array}$$
%If $S$ is not isomorphic to a direct product with a cyclic group, it will thus prove the following result.
%\begin{Theorem}
%The Freedman group $Fr(162\times 4)$ contains a subgroup which is isomorphic to $\Delta(54)$.
%\end{Theorem}
%\newtheorem{Remark}{Remark}
%\begin{Remark}
%In $2009$, R. Zwicky and T. Fischbacher showed in their paper \cite{FZ} that every $(D)$-group is a subgroup of $\Delta(6\,n^2)$ for a suitable $n$. Then the Freedman group also arises as a subgroup of some $\Delta(6n^2)$.
%\end{Remark}
In what follows, decide to set $$S=\,(<A>\times <B>)\,\rtimes <H_2>$$
By the same arguments as in the proof of Lemma $3$, all the elements of odd order of $S$ belong to $\mathcal{N}$.
To proceed, we will need another fact, this time about the elements of order $2$.
\begin{Lemma}
There are exactly three elements of order $2$ in $S$, namely $H_2$, $BH_2$ and $B^2H_2$.
\end{Lemma}
\begin{Corollary}
The group $S$ is isomorphic to the direct product $\mathbb{Z}_9\times S_3$.
\end{Corollary}
\textsc{Proof.} (Corollary) In the dihedral group $D_9$, there are exactly $9$ elements of order $2$. Thus, the number of elements of order $2$ in $S$ is not sufficient for the group to be isomorphic to $\mathbb{Z}_3\times D_9$. \\

\textsc{Proof.} (Lemma) We know from studying the structure of the Freedman group that the elements of $S$ can be uniquely written as products $A^iB^jH_2$ with $i\in[\!|0,8|\!]$ and $j\in\lbrace 0,1,2\rbrace$. So, an element has order two means
\begin{equation}H_2A^iB^jH_2=A^{-i}B^{-j}\end{equation}
From there, since
\begin{eqnarray}
H_2A\,H_2^{-1}&=&AB\\
H_2\,B\,H_2^{-1}&=&B^2
\end{eqnarray}
we claim that the values $i=1,2,4,5,7,8$ are to exclude to allow Eq. $(57)$ to be verified. Indeed, we have the following lemma.
\begin{Lemma}
Out of the $18$ elements of order $9$
$$\begin{array}{l}
A,\,A^2,\,A^4,\,A^5,\,A^7,\,A^8\\
AB,\,A^2B,\,A^4B,\,A^5B,\,A^7B,\,A^8B\\
AB^2,\,A^2B^2,\,A^4B^2,\,A^5B^2,\,A^7B^2,\,A^8B^2
\end{array}$$
of $S$, the ones that are self-conjugate under $H_2$ are exactly those below.
$$AB^2,A^2B,A^4B^2,A^5B,A^7B^2,A^8B$$
These are the only elements of $S$ which have order $9$ and belong to $Z(S)$.
\end{Lemma}
\noindent This Lemma follows from the key relation
\begin{equation}
H_2\,A^iB^j\,H_2=A^iB^{i+\bar{j}}
\end{equation}
which holds for any $i$ and $j$ with $i\in[\!|0,8|\!]$ and $j\in\lbrace 0,1,2\rbrace$. In this equality, $\bar{j}$ denotes the "conjugate" of $j$, that is $\bar{j}=2$ if $j=1$ and $\bar{j}=1$ if $j=2$.\\
Indeed, it suffices to notice that
$$\forall\,(i,j)\in\lbrace (1,2),(2,1),(4,2),(5,1),(7,2),(8,1)\rbrace,\;i+\bar{j}\;\equiv\; j\;(mod\,3)$$
\hfill $\square$\\
We now return to the proof of Lemma $14$. \\
By Lemma $15$, for the values of $i$ and $j$ above, the product matrix $A^iB^j$ is self-conjugate under $H_2$. Moreover, it cannot equal its inverse since it has order $9$. Thus, Eq. $(57)$ is not satisfied for such $i$'s and $j$'s. And still for the same values of $i$, Eq. $(57)$ is still not satisfied when $j=0$. This follows from Eq. $(58)$ and the fact that for these given values of $i$, we have $B^i\in\lbrace B,B^2\rbrace$. Now the contradiction comes from
\begin{equation}<A>\cap<B>\,=\,\lbrace I_3\rbrace\end{equation}
Also, for these values of $i$ and the respective conjugate values for $j$, the $H_2$ conjugate of $A^iB^j$ is $A^i$. Again, Eq. $(61)$ prevents Eq. $(57)$ from happening.
Thus, it remains to check the $H_2$-conjugates of the elements of $\mathcal{N}$ of order $3$. From Eq. $(60)$, these are ruled by the relation
\begin{equation}H_2A^iB^jH_2=A^iB^{\bar{j}}\,,\end{equation}
as now $i\in\lbrace 0,3,6\rbrace$.
%where $\bar{j}$ denotes the conjugate of $j$ like previously defined. \\
By $(61)$, both Eqs. $(57)$ and $(62)$ imply that $i$ must then be zero. And since $j+\bar{j}=3$, Eq. $(57)$ gets
satisfied when $i=0$ and $j\in\lbrace 1,2\rbrace$, thus ending the proof of Lemma $14$.
\begin{Proposition}
$$S=\,<AB^2>\times(<B>\rtimes<H_2>)\simeq\mathbb{Z}_9\times S_3$$
\end{Proposition}
\begin{Corollary}
There exists semi-direct products so that
$$(\mathbb{Z}_9\times\mathbb{Z}_3)\rtimes\mathbb{Z}_2\simeq\mathbb{Z}_9\times(\mathbb{Z}_3\rtimes\mathbb{Z}_2)$$
\end{Corollary}
\textsc{Proof.} By Lemma $14$, the group $S$ contains a subgroup $<B>\rtimes<H_2>$ which is isomorphic to $S_3$. Further, by Lemma $15$, the product $AB^2$ is self-conjugate under $H_2$ and thus commutes to all the other elements in the group. Moreover, we have
$$S=\,<AB^2>.<B>.<H_2>$$
Then $S$ is the direct product announced in the statement of Proposition $6$. \\
Corollary $6$ follows. \hfill $\square$\\

We now deal with the proof of Theorem $13$. Again, it relies on the table from the link \cite{PW2} in \cite{PW}. If $Fr(162\times 4)$ were isomorphic to a finite subgroup of $U(2)$, then the four isomorphic $3$-Sylow subgroups of $Fr(162\times 4)$ of order $81$ would also be isomorphic to a finite subgroup of $U(2)$. However, when looking at the short list of non-abelian groups of order $81$ in the table, it appears that the only non-abelian $SU(3)$ finite subgroup of that order is
$$C(9,1,1)=(\mathbb{Z}_9\times\mathbb{Z}_3)\rtimes\mathbb{Z}_3,$$
which, as read from the table, cannot be identified with a finite subgroup of $U(2)$. \hfill $\square$

\section{Conclusion}
\subsection{Computer verification with GAP}
Around the $2000$ millenium time, Hans Besche, Bettina Eick and Eamonn O'Brien announced the construction up to isomorphism of the $49\, 910\,529\,484$ groups of order at most $2000$. Their work is of course amazing and from their Table $1$ of \cite{OB}, the most difficult orders (in terms of the largest numbers of groups) appear as some products of a power of $2$ by a power of $3$.  We know from their work \cite{OB} that there are $757$ non-isomorphic groups of order $648$. We provide below a small table which contains a single row of their big table. \\\\
\begin{tabular}{|c|c|c|c|c|c|c|c|c|c|c|}
\hline
\text{Order}&640&641&642&643&644&645&646&647&648&649\\
\hline
$\begin{array}{l}\text{Number}\\\text{of groups}\end{array}$&21541&1&4&1&9&2&4&1&757&1\\
\hline
\end{tabular}
$$\begin{array}{l}\end{array}$$
We entered the presentation in GAP and found out that $Fr(162\times 4)$ is the group of GAP ID \begin{equation}[648,259]\end{equation} As for the structure returned by GAP for this group ID, it shows
\begin{equation}((\mathbb{Z}_{18}\times\mathbb{Z}_{6})\rtimes\mathbb{Z}_3)\rtimes\mathbb{Z}_2\end{equation}
%And so our study shows that there exists an isomorphism
%\begin{equation}(\mathbb{Z}_{18}\times\mathbb{Z}_6)\rtimes\mathbb{Z}_3)\rtimes\mathbb{Z}_2\simeq (\mathbb{Z}_{18}\times\mathbb{Z}_6)\rtimes(\mathbb{Z}_3\rtimes\mathbb{Z}_2)\end{equation}
We then entered the presentation determined in the present paper for the group $D(18,1,1;2,1,1)$ with the reassuring outcome that both groups have the same ID, thus confirming that they are isomorphic.
%We last entered the presentation read out of \cite{GL} for the exceptional group $\Sigma(216\times 3)$ and found a different GAP ID, namely $3$
We also played the same game with the smaller groups $D(9,1,1;2,1,1)$ and $Fr(162)$. A presentation for the first group is part of the material exposed in the present paper, and a presentation for the second group can be found at the end of \cite{BL}. Note that one of the relations in the presentation of \cite{BL} is redundant and the commutator relation was unwillingly forgotten by the authors. Both groups have the same GAP ID, that is
\begin{equation}
[162,14]
\end{equation}
Notice -- but only for the amusement and the beauty -- that
$$\begin{array}{l}\text{$Fr(162)$ is the $7\times 2$-th group of order $162$}\\
\text{$Fr(162\times 4)$ is the $7\times 37$-th group of order $648$}\end{array}$$
 out of $55$ groups of order $162$ for the first one and out of $757$ groups of order $648$ for its extension.
% (out of $55$ groups of order $162$) and $Fr(162\times 4)$ is the $7\times 37$-th group of order $648$ (out of $757$ groups of order $648$).
Note that there is a composite number of groups of order $162$ and a prime number of groups of order $648$.
For this ID, GAP provides the structure \begin{equation}((\mathbb{Z}_9\times\mathbb{Z}_3)\rtimes\mathbb{Z}_3)\rtimes\mathbb{Z}_2\end{equation}
\newtheorem{Remark}{Remark}
\begin{Remark}
Combining the results from GAP and our own study, we see that we obtain some isomorphisms:
$$\begin{array}{cc}
((\mathbb{Z}_{18}\times\mathbb{Z}_6)\rtimes\mathbb{Z}_3)\rtimes\mathbb{Z}_2\simeq(\mathbb{Z}_{18}\times\mathbb{Z}_6)
\rtimes(\mathbb{Z}_3\rtimes\mathbb{Z}_2)&\text{by Eq. $(64)$ and $\S\,2$}\\
&\\
((\mathbb{Z}_9\times\mathbb{Z}_3)\rtimes\mathbb{Z}_3)\rtimes\mathbb{Z}_2\simeq(\mathbb{Z}_9\times\mathbb{Z}_3)\rtimes
(\mathbb{Z}_3\rtimes\mathbb{Z}_2)&\text{by Eq. $(66)$, \cite{BL} and $\S\,3$}\end{array}$$
\end{Remark}
\noindent Associativity of the semi-direct product up to isomorphism is not automatic. But holds for instance when $G_1$, $G_2$, $G_3$ are subgroups of a same group $G$ and the group $G_2\rtimes G_3$ acts on $G_1$ by first conjugating by the element of $G_3$ and then conjugating by the element of $G_2$. More generally, if
$$\begin{array}{cccccccc}
\varphi:&G_3&\rightarrow&Aut(G_1\rtimes G_2),&\psi:&G_2&\rightarrow&Aut(G_1),\\
\alpha:&G_3&\rightarrow&Aut(G_2),&\beta:&G_2\rtimes G_3&\rightarrow& Aut(G_1),
\end{array}$$
are homomorphisms satisfying to
\begin{eqnarray*}
\forall\,(g_2,g_3)\in G_2\times G_3,\qquad\;\;\beta(g_2,g_3)\circ p_1&=&\psi(g_2)\circ p_1\circ\varphi(g_3)\\
\forall\,(g_2,g_3)\in G_2\times G_3,\qquad\qquad\alpha(g_3)\circ p_2&=&p_2\circ\varphi(g_3)
\end{eqnarray*}
where $p_1$ and $p_2$ denote the respective projections with respect to the first and second coordinates of the cartesian product, then we have
$$(G_1\underset{\psi}{\rtimes}G_2)\underset{\varphi}{\rtimes}G_3\simeq G_1\underset{\beta}{\rtimes}(G_2\underset{\alpha}{\rtimes}G_3)$$
%which agrees with our own "hand-made work" to determine the structure of the group since it is easy %proves the existence of an isomorphism
%\begin{equation}
%((\mathbb{Z}_9\times\mathbb{Z}_3)\rtimes\mathbb{Z}_3)\rtimes\mathbb{Z}_2\simeq(\mathbb{Z}_9\times\mathbb{Z}_3)\rtimes (\mathbb{Z}_3\rtimes\mathbb{Z}_2)
%\end{equation}
%\begin{center}
%To be completed.
%\end{center}
\subsection{Concluding remarks}
One of the weaknesses of the classification of \cite{MB} or its extended version is that this classification is not a classification up to isomorphism as well pointed out during our mini studies involving the $SU(3)$ finite subgroups from the series $(D)$. Also,
a systematic analysis of the structure of some of the groups is not yet complete. The classification of all the finite subgroups of $SU(3)$ up to isomorphism remains an open problem. However, a stronger classification, this time up to conjugacy, makes more sense when dealing with groups of matrices. Of course a necessary condition for conjugacy is isomorphism.
We investigated whether the Freedman group $Fr(162\times 4)$ is conjugate to $D(18,1,1;2,1,1)$ and found out this answer is yes. The conclusions we drew are summarized below.
\begin{Theorem}\hfill\\
$(i)$ The Freedman group $Fr(162\times 4)$ is conjugate to $D(18,1,1;2,1,1)$.\\
$(ii)$ There exists an orthogonal matrix $O$ such that
$$O\,Fr(162\times 4)O^{T}=D(18,1,1;2,1,1)$$
with $$O=\begin{pmatrix} 1/\sqrt{2}&0&1/\sqrt{2}\\0&1&0\\-1/\sqrt{2}&0&1/\sqrt{2}\end{pmatrix}\\$$
$(iii)$ The same conjugation relation holds for the respective subgroups: $$O\,Fr(162)\,O^{T}=D(9,1,1;2,1,1)$$
\end{Theorem}
\textsc{Proof.} Suppose that there exists an invertible matrix $P\in GL_3(\mathbb{C})$ such that $$P\,[Fr(162\times 4)]P^{-1}=D(18,1,1;2,1,1)$$
\begin{Lemma}
The invertible matrix $P$ must be a transition matrix from a common basis of diagonalization for $\mathcal{N}$ to the canonical basis of $\mathbb{C}^3$.
\end{Lemma}
\textsc{Proof.} Under conjugation, a $3$-Sylow subgroup of $Fr(162\times 4)$ must be mapped to a $3$-Sylow subgroup of $D(18,1,1;2,1,1)$. In particular, we must have
\begin{equation}
P\mathcal{S}_3(F)P^{-1}=W_i\mathcal{S}_3(D)W_i^{-1}\;\;\text{some}\;i\in\lbrace 1,2,3,4\rbrace\;\text{with $W_1=I_3$}
\end{equation}
So, we have
\begin{equation}
P\mathcal{S}_3(F)P^{-1}=\mathcal{F}\sqcup\mathcal{F}W_iEW_i^{-1}\sqcup\mathcal{F}W_iE^2W_i^{-1}
\end{equation}
The key idea now is to notice that the traces of elements of $\mathcal{F}W_iEW_i^{-1}$ or $\mathcal{F}W_iE^2W_i^{-1}$ are all zero, this for all $i\in\lbrace 1,2,3,4\rbrace$. Since the trace of $A$ is not zero, we must then have
\begin{equation}
PAP^{-1}\in\mathcal{F}
\end{equation}
The same trick does not work with $B$ because $B$ has trace zero. However, the trace of the product $AB$ is non-zero.
Therefore,
\begin{equation}
P(AB)P^{-1}=(PAP^{-1})(PBP^{-1})\in\mathcal{F}
\end{equation}
Now $(69)$ and $(70)$ imply that $PBP^{-1}$ is also in $\mathcal{F}$. Hence the immediate Corollary.
\begin{Corollary}
\begin{equation}
P\mathcal{N}P^{-1}=\mathcal{F}
\end{equation}
\end{Corollary}
\noindent This finishes the proof of Lemma $16$. \hfill$\square$ \\\\
%And since all the matrices of $\mathcal{N}$ commute, we see that $P$ is a transition like announced in the statement of Lemma.
\noindent To move further, we will need to use some facts from before. We have seen in the proof of Theorem $7$ that any isomorphism of groups between $Fr(162\times 4)$ and $D(18,1,1;2,1,1)$ must map the Klein group $\mathcal{V}$ onto the Klein group $\mathcal{V}_D$.
We apply this fact to our isomorphism by conjugation and we get
\begin{equation}
P\mathcal{V}P^{-1}=\mathcal{V}_D
\end{equation}
It now suffices to recall that
$$\begin{array}{ccccc}
(\mathcal{N}\times\mathcal{V})&\cap &S_3(F)&=&\lbrace I_3\rbrace\\
(\mathcal{F}\times\mathcal{V}_D)&\cap&S_3(D)&=&\lbrace I_3\rbrace
\end{array}$$
Then, we also have
\begin{equation}
PS_3(F)P^{-1}=S_3(D)
\end{equation}
In order to conclude, we can't bypass to investigate the shape of the transition matrix $P$. The matrix $A$ has three distinct eigenvalues, hence each eigenspace has dimension $1$ over $\mathbb{C}$. Explicitly, the eigenspaces are
$$\mathbb{C}\,\left(\begin{array}{l}-1\\0\\1\end{array}\right),\,\mathbb{C}\,\left(\begin{array}{l}0\\1\\0\end{array}
\right),\,\mathbb{C}\,\left(\begin{array}{l}1\\0\\1\end{array}\right)$$
Then, the most general form for $P$ is as follows.
\begin{equation}
P=\left[\begin{array}{l}\begin{pmatrix} -1&0&1\\0&1&0\\1&0&1\end{pmatrix}\begin{pmatrix}d_1&&\\&d_2&\\&&d_3\end{pmatrix}P_{\sigma}\end{array}\right]^{-1}\,,
\end{equation}
where $d_1$, $d_2$, $d_3$ are non-zero complex numbers and $P_{\sigma}$ is the permutation matrix associated with a given permutation $\sigma$ of $Sym(3)$. \\
Equality $(73)$ imposes some restrictions on $d_1$, $d_2$, $d_3$ and $\sigma$.
% For instance, the matrix
%$P(FUM)^3P^{-1}$ must be diagonal and moreover one of the $W_i$'s with $i\in\lbrace 1,2,3,4\rbrace$.
First, we computed the product
\begin{equation}
\begin{pmatrix}1/d_1&&\\&1/d_2&\\&&1/d_3\end{pmatrix}\begin{pmatrix}-1&0&1\\0&1&0\\1&0&1\end{pmatrix}^{-1}H_3
\begin{pmatrix}-1&0&1\\0&1&0\\1&0&1\end{pmatrix}\begin{pmatrix}d_1&&\\&d_2&\\&&d_3\end{pmatrix}
\end{equation}
and found out this product is:
$$\begin{pmatrix} 0&-\frac{d_2}{\sqrt{2}\,d_1}&0\\
0&0&\frac{\sqrt{2}d_3}{d_2}\\
-\frac{d_1}{d_3}&0&0
\end{pmatrix}$$
This matrix must also equal one of $E$ or $\overset{\sim}{B}E\overset{\sim}{B}$. By looking at the shape of all these matrices and also allowing conjugating Eq. $(75)$ by a permutation matrix, we see that the following conditions must hold anyway.
$$\left\lbrace\begin{array}{l}
d_3=-d_1\\
d_2=-\sqrt{2}d_1
\end{array}\right.$$
That is all the non-zero coefficients must be $1$'s. Conversely, \begin{center}
set $\left|\begin{array}{ccc}d_1&=&-\frac{1}{\sqrt{2}}\\d_2&=&1\\d_3&=&\frac{1}{\sqrt{2}}\\\sigma&=&id\end{array}\right.$
\end{center}
Then, we have
\begin{eqnarray*}
PH_3P^{-1}&=&E\\
PH_1P^{-1}&=&\overset{\sim}{B}E\\
P(FUM)^3P^{-1}&=&W_3\\
PG_1(FUM)^3G_1^{-1}P^{-1}&=&W_4\\
P\mathcal{N}P^{-1}&=&\mathcal{F}
\end{eqnarray*}
so that the first two points of Theorem $14$ hold. As for point $(iii)$, it is part of our latter work since for this $O$, we have just seen that
$$\left\lbrace\begin{array}{ccccc}
O&\mathcal{N}&O^{T}&=&\mathcal{F}\\
O&S_3(F)&O^{T}&=&S_3(D)
\end{array}\right.$$
%We still want to justify that $Fr(162\times 4)$ is conjugate to neither of $D(18,1,1;2,0,0)$ and $D(18,1,1;2,0,1)$. For the first one, an easy argument is to notice that both groups differ only by a minus sign in front of the generator $\overset{\sim}{B}$. So, conjugating by $\overset{\sim}{B}$ or by $-\overset{\sim}{B}$ is equivalent. Then, whatever be the choice of permutation for $\sigma$, the same relations as before must hold on the $d_i$'s and forbid to have
%$$PH_1P^{-1}\in\lbrace -\overset{\sim}{B},-\overset{\sim}{B}E,-E\overset{\sim}{B}\rbrace$$
%As for $D(18,1,1;2,0,1)$,
%and
%$$\begin{array}{l}
%\text{$PAP^{-1}$ is diagonal}\\
%\text{$PBP^{-1}$ is diagonal}
%\end{array}$$
%Conjugating a diagonal matrix by any permutation matrix has the effect of doing the same operation on the rows and on the columns which results in permuting the coefficients on the diagonal.
\section{Appendix}
This part presents a brief explanation for the non-expert for how to obtain the permutation matrix from the Freedman fusion operation. We use the unitary normalization of Kauffman-Lins theory at level $4$ and refer the unfamiliar reader to \cite{KL} and \cite{ZW}. Here is how we evaluate the diagram
\begin{center}
\epsfig{file=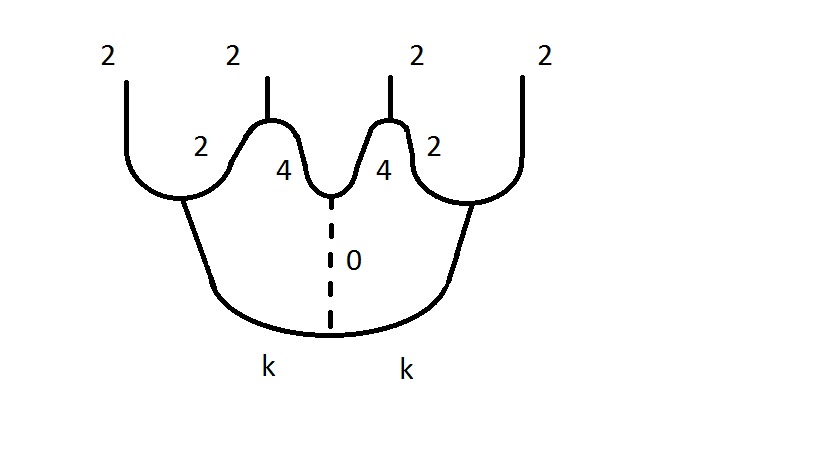, height=4.5cm}
\end{center}
First, we do an $F$-move on the edge labeled "0". From the fusion rules, set
$$\begin{array}{ccc}
i=4&\text{when}&k=0\\
i=2&\text{when}&k=2\\
i=0&\text{when}&k=4
\end{array}$$
We obtain\\
\epsfig{file=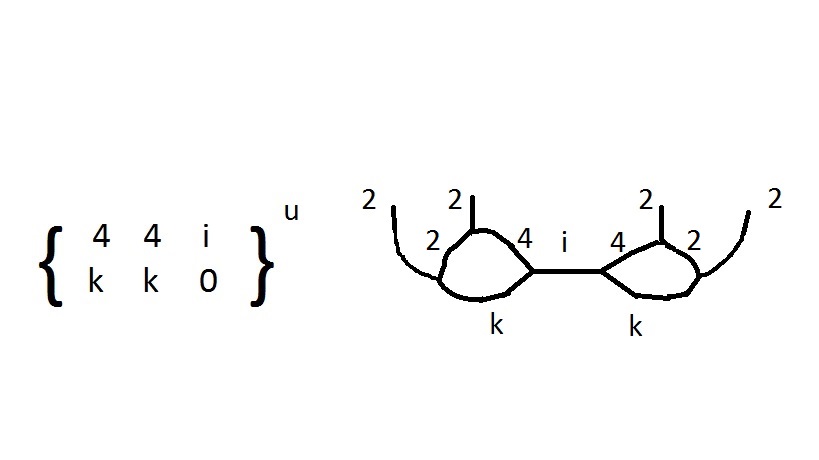, height=6cm}\\
where the brackets are used to denote unitary $6j$-symbols. Do two more $F$-moves on the "external" edges labeled "2". Get
\begin{center}
\epsfig{file=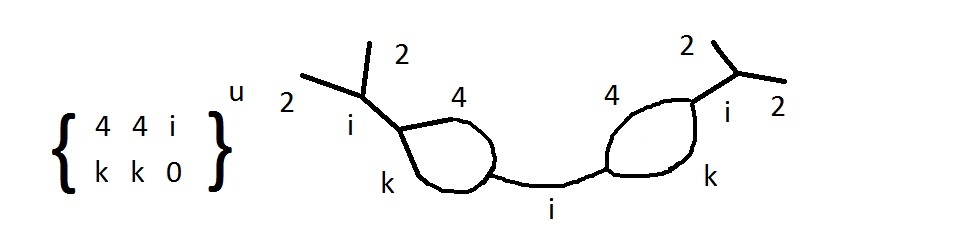, height=3.5cm}
\end{center}
Next, undo the two loops by multiplying by adequate unitary theta symbols divided by the quantum dimension $\Delta_i$ of the particle of topological charge $i$. Get
\begin{center}\epsfig{file=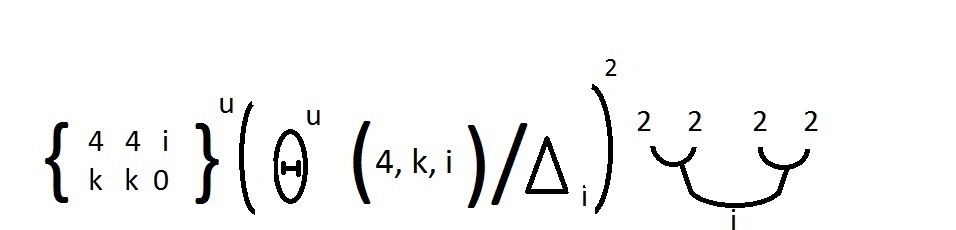, height=2.5cm}\end{center}
$\begin{array}{l}\\\end{array}$ \\
We see that the particles of respective topological charge $0$ and $4$ have been interchanged as a result of the Freedman fusion operation. Recall from \cite{ZW} that
$$\Theta^u(a,b,c)=\sqrt{\Delta_a}\sqrt{\Delta_b}\sqrt{\Delta_c}$$
Also, the unitary version of the $6j$-symbol is the following.
$$\left\lbrace\begin{array}{ccc} G&B&E\\C&D&F\end{array}\right\rbrace^u=\frac{Tet\left[\begin{array}{ccc} G&B&E\\C&D&F\end{array}\right]\sqrt{\Delta_E}\sqrt{\Delta_F}}{\sqrt{\Theta(G,D,E)}\sqrt{\Theta(C,D,F)}\sqrt{\Theta(C,B,E)}
\sqrt{\Theta(G,B,F)}}$$
where $Tet$ is the tetrahedron of \cite{KL} and $\Theta$ denotes the non-unitary theta symbol. \\\\
We evaluated the coefficient in front of the final diagram and found the value $1$ for all the couples
$(k,i)=(0,4)$, $(k,i)=(2,2)$ and $(k,i)=(4,0)$. Thus, the resulting unitary matrix of the Freedman fusion operation  is
$$\begin{array}{cc}\begin{array}{l}\\\\ |0>\\|2>\\|4>\end{array}&\begin{array}{l}
\qquad\negthickspace|0>\, |2>\, |4>\\\\
\begin{pmatrix}& 0&&0&&1&\\
&0&&1&&0&\\
&1&&0&&0&\end{pmatrix}
\end{array}\end{array}$$
As a matrix of $SU(3)$ we obtain the fusion matrix which we named $FUM$.\\\\
\begin{center}\textbf{Acknowledgements}\end{center}
The author thanks Michael Freedman for his generosity and enthusiasm at communicating his ideas, thus allowing this paper to grow. She thanks Zhenghan Wang for introducing her to the Kauffmann-Lins theory with much kindness and generosity. Thanks go to Bela Bauer for numerically confirming the order of the Freedman group. The author heartfully thanks Patrick Otto Ludl for getting back to her with much kindness and at the speed of light when she several times seeked his expertise on the finite subgroups of $SU(3)$. She is pleased to thank Eamonn O'Brien for informative communications about the possibilities offered by GAP. She thanks Michael Freedman for helpful discussions and David Wales for wise comments. Thanks go to Jennifer Cano, Meng Cheng, Sean Fraer, Younghyun Kim and Eugeniu Plamadeala for their friendship during the preparation of this manuscript and the whole team of Microsoft Research Station Q for offering such a pleasant environment of work.

\end{document}